\magnification=1000
\input amstex
\input pictex
\documentstyle{amsppt}
\loadbold \baselineskip = 18pt \overfullrule = 0pt

\vsize=9.0truein \hsize=6.5truein \voffset=-.5in

\baselineskip = 18pt


\def\cl#1{\Cal{#1}}

\def\ovl{\overline}
\def\ds{\displaystyle}


\def\dc#1{\partial C(L^#1)}

\def\siml#1{{\overset{L^{\#}}\to\sim}}

\def\re{\operatorname{Re}}

\def\im{\operatorname{Im}}

\def\ad{\operatorname{ad}}

\def\tr{\operatorname{tr}}

\def\det{\operatorname{det}}
\def\diag{\operatorname{diag}}

\def\norm#1{\|#1\|}

\def\pmtwo#1#2#3#4{\pmatrix#1&#2\\#3&#4\endpmatrix}

\def\cb{{\Bbb C}}
\def\rb{{\Bbb R}}

\let\loz=\lozenger

\let\tildesymbol=\~

\def\cht{C_{\hat w_\theta}}
\def\ct{C_{w_\theta}}

\def\cxs{C_{w_x}}

\def\epts{e^{i\theta\ads}}

\let\ads=\adsg

\def\norm#1#2{\|#2\|_{L^{#1}}}
\def\Norm#1#2#3{\|#3\|_{L^{#1}(#2)}}

\def\ep{\epsilon}

\let\loz=\lozenge

\def\ads{\ad\sigma\,}

\def\sr{{\Cal R}}

\define\vp{\epsilon}
\define\n{\noindent}

\NoBlackBoxes
\leftheadtext{P. Deift and X. Zhou} \rightheadtext{Long--time
asymptotics }

\topmatter\title Long--time asymptotics for solutions of the NLS
equation\\ with initial data in a weighted Sobolev
space\endtitle\endtopmatter

\document
\centerline{Percy Deift} \centerline{New York University}
\vskip1cm \centerline{ Xin Zhou} \centerline {Duke University}

\vskip2cm

\head \S 1. Introduction\endhead

Long--time asymptotics for solutions $q(x,t)$ of the defocusing
Nonlinear Schr\"odinger (NLS) Equation
$$\left\{\aligned &iq_t + q_{xx} - 2 |q|^2q = 0\\
&q(t=0,x) = q_0(x) \to 0 \text{ as } |x|\to \infty
\endaligned\right. \tag1.1$$
\noindent have been obtained in [ZaMa][DIZ][DZ2][DZ4] for initial
data with sufficient smoothness and decay: as $t \to \infty$
$$q(x,t)=t^{-1/2}\alpha(z_0)e^{ix^2/(4t)-i\nu(z_0)\log2t} +O(\log
t/t),\tag1.2$$where
$$\nu(z)=-\frac1{2\pi}\log(1-|r(z)|^2),\ \ |\alpha(z)|^2=\nu(z)/2,$$and
$$\arg\alpha(z)=\frac1\pi\int_{-\infty}^{z}\log(z-s)d(\log(1-|r(s)|^2))+\frac\pi4
+\arg\Gamma(i\nu(z))+\arg r(z).$$ Here $\Gamma$ is the gamma
function and the function $r$ is the so-called reflection
coefficient for the potential $q_0(x)=q(x,t=0)$, as described
below. The error term $O(\frac{\log t}t)$ is uniform for all
$x\in\rb$.  The above asymptotic form was first obtained in
[ZaMa], but without the error estimate. Based on the nonlinear
steepest descent method introduced in [DZ1], the error estimate in
(1.1) was derived in [DIZ] (see also [DZ2] for a pedagogic
presentation). As noted above, some high orders of decay and smoothness are
required for the initial data.  In this paper we describe a new
method that produces an error estimate of order
$O(t^{-\frac12-\kappa})$ for any $0<\kappa<\frac14$, just under
the assumption that the initial data $q_0$ lies in the weighted
Sobolev space $H^{1,1}=\{f\in L^2(\rb):\ xf,f'\in L^2(\rb)\}$.
Such an estimate is needed, for example, in [DZ3][DZ5] where the
authors obtain long-time asymptotics for solutions of the
perturbed NLS equation, $iq_t + q_{xx} - 2|q|^2q - \vp |q|^lq=0$,
for $l>2$ and $\vp>0$. As we will see (cf. Section 4 below), the
estimate $O(t^{-\frac12-\kappa})$ in fact depends only on the
weighted $L^2$ norm of the initial data,
$\left(\int_\rb(1+x^2)|q_0(x)|^2dx\right)^{1/2}<\infty.$ The
estimate is (essentially) optimal even in the linear case where
the standard Fourier method produces an error of order
$O(t^{-3/4})$. Our new method is a further development of the
steepest descent method of [DZ1], and replaces certain key
absolute type estimates in [DZ1] with cancellations from
oscillations.

The NLS equation can be integrated by using the familiar
scattering theory/inverse scattering theory for the ZS--AKNS
system [ZS][AKNS] associated to NLS,
 $$\partial_x \psi =
\left(iz\sigma + \left(\matrix 0&q \\ \bar q&0\endmatrix
\right)\right)\psi\tag1.3$$
where $\sigma = \left( \matrix 1/2&0\\
0&-1/2\endmatrix \right)$. As is well known [ZS], the NLS equation
is equivalent to an isospectral deformation of the operator
$\partial_x - \left(iz\sigma + \left(\matrix 0&q\\ \bar
q&0\endmatrix\right)\right)$. As described in Section~3 below, for
each $z\in {\Bbb C}\backslash {\Bbb R}$, one constructs solutions
$\psi(x,z)$ of (1.3) of the type considered in [BC]  with the
following properties:

\itemitem{(a)} $m(x,z) \equiv
\psi(x,z) e^{-ixz\sigma} \to I = \left(\matrix 1&0\\
0&1\endmatrix\right)$ as $x\to -\infty$,
\itemitem{(b)} $m(x,z)$ is bounded as $x\to +\infty$.

\n For each fixed $x$, the $2\times 2$ matrix function $m(x,z)$
solves the following Riemann-Hilbert problem (RHP) in $z$:
\flushpar (1.4) \ (i) \ $m(x,z)$ is analytic in  ${\Bbb C}
\setminus {\Bbb R}$,

\ \ \ \ \ \ (ii) \  $m_+(x,z) = m_-(x,z) v_x(z),\qquad z\in {\Bbb
R}$, \flushpar where $ m_\pm (x,z) = \lim_{\vp\downarrow 0} m(x,z
\pm i\vp)$ and $v_x(z) = \left(\matrix 1-|r(z)|^2&r(z)e^{izx}\\
-\overline{r(z)}e^{-izx}&1\endmatrix\right)$ where $r(z)$ is the
{\it reflection coefficient} of $q$.

\ \ \ \ \ \ (iii) \  $\lim\limits_{z\to \infty} m(x,z) = I$.

\n The sense in which the limits in (ii) and (iii) are achieved
will be made precise in Sections~2 and 3. The reflection
coefficient satisfies the important {\it a priori\/} bound
$\|r\|_{L^\infty(dz)} <1$.

If we expand out the limit in (iii),
$$m(x,z) = I+ \frac{m_1(x)}z + o\left(\frac1{z} \right),\tag1.5$$
then we obtain an expression for $q$,
$$q(x) = -i(m_1(x))_{12}.\tag1.6$$
The direct scattering map ${\Cal R}$ is obtained by mapping
$q\mapsto r$,
$$q\mapsto m(x,z) = m(x,z;q) \mapsto v_x(z) \mapsto r =  {\Cal R}(q).$$
Given $r$, the inverse scattering map ${\Cal R}^{-1}$ is obtained
by solving the RHP (1.4) and mapping to $q$ via (1.6),
$$r\mapsto \text{ RHP } \mapsto m(x,z)  = m(x,z;r) \mapsto m_1(x) \mapsto q = {\Cal R}^{-1}(r).$$

The remarkable fact discovered in [ZS] is that if $q(t)=q(x,t)$
solves the NLS equation, then $r(t)={\Cal R}(q(t))$ evolves simply
as $$ r(t)=r(z,t)=e^{-itz^2}r_0(z),\tag1.7 $$ where $r_0={\Cal
R}(q_0=q(t=0))$. Using $\loz$ to denote a dummy variable, we may
rewrite (1.7) as a formula $$ q(t)={\Cal
R}^{-1}(e^{-i\loz^2t}r(\loz,t=0)),\tag1.8$$ for the solution of
NLS. This formula shows that the problem of computing the
asymptotics of $q(t)$ as $t\to \infty$, reduces to the problem of
analyzing the map $\sr^{-1}$, and it was precisely for this
purpose that the nonlinear steepest descent method was introduced
in [DZ1] in the context of the modified Korteweg deVries equation.

In order to proceed, we must specify the domain and range of
$\sr$. Many authors have commented that $\sr$ is a nonlinear
Fourier-type map and indeed the results in [BC] imply that $\sr$
is a bijection from Schwartz space $\Cal S(\rb)$ onto ${\Cal
S}_1(\rb))={\Cal S}(\rb)\cap\{r:\ \|r\|_{L^\infty(\rb)}<1\}.$ The
first proof that $\sr,\ \sr^{-1}$ involve no ``loss" of smoothness
or decay, was given in [Z1] where the author showed that $\sr$ is
a bijection (in fact a bi-Lipschitz mapping) from the weighted
Sobolev space $$ H^{k,j}=\{f:\ f,\partial_x^kf,x^jf \in
L^2(\rb)\}$$ with norm
$$\|f\|_{H^{k,j}} = (\|f\|^2_{L^2} + \|\partial^k_x
f\|^2_{L^2} +
 \|x^jf\|^2_{L^2})^{\frac12}.$$
onto $H_1^{j,k}=H^{j,k}\cap\{f:\|f\|_{L^\infty(\rb)}<1\}$  for
$k\ge0,j\ge1$.

We say that $q(t),\ t\ge0$, is a {\it (global) weak solution} of
(1.1) in $H^{k,j}$ if $t \mapsto q(t)$ is a continuous map from
$\rb_+$ into $H^{k,j}$ satisfying$$q(t) = e^{-iH_0t} q_0 - i
\int^t_0 e^{-iH_0(t-s)} 2|q(s)|^2 q(s)ds\tag1.9$$ for $t\ge0$,
where $H_0=-\partial_x^2$ is negative Lapalcian regarded as a
self-adjoint operator on $L^2(\rb)$. Standard estimates (see, for
example, [DZ5] and the references therein) show that weak
solutions of (1.1) exist and are unique in $H^{k,j}$ for all
$k\ge1,j\ge0$. Thus the largest $H^{k,j}$ space for which weak
solutions $q(t)$ exist and for which $r(t)=\sr(q(t))$ is defined
is $H^{1,1}$. Henceforth, by a solution $q(t)$ of (1.1) we always
mean the (unique, global) weak solution in $H^{1,1}$. By the
preceding remarks, $r(t)=\sr(q(t))$ is given by
$e^{-iz^2t}r(z,t=0)$ and we verify directly that if
$r(t=0)=\sr(q_0)\in H_1^{1,1}$, then so does $r(t)$ for all
$t\ge0$.

Our main result is the following (see also Theorem 4.40 below).

\proclaim{Theorem 1.10}Let $q(t),\ t\ge0$, solve (1.1) with
$q_0=q(t=0)\in H^{1,1}$. Fix $0<\kappa<1/4$. Then as
$t\to\infty$,$$q(x,t)=t^{-1/2}\alpha(z_0)e^{ix^2/(4t)-i\nu(z_0)\log2t}+
O\left(t^{-(1/2+\kappa)}\right),$$ where $\alpha$ and $\nu$ are
given in terms of $r=\sr(q_0)$ as above. The error term $
O\left(t^{-(1/2+\kappa)}\right)$ is uniform for all $x\in\rb$.
\qed\endproclaim

The paper is organized as follows. In \S2 we make the notion of
the RHP in (1.4) precise, and also describe some associated
inhomogenous RHP's which are useful in proving Theorem 1.10. In
particular, we use the ideas of this section to prove that the
resolvent of the solution operator for an associated model RHP is
uniformly bounded in space and time as $t\to\infty$ (see
Proposition 2.49 below). Some parts of this section are discussed
on the WEBPAGE of [DZ5]. In \S3 we give a proof that $\sr$ is a
bi-Lipschitz map from $H^{1,1}$ onto $H^{1,1}_1$. As in the case
of the long-time behavior of $q(t)$ noted above, the bijectivity
of $\sr$ relies on cancellations from oscillations, rather than
brute force estimation. Our proof follows [Z1] where the author
considers a very general class of systems for all (allowable)
$k,j$. For such general systems many additional technicalities
arise, and for this reason we believe it is useful for the reader
for us to present in detail the bijectivity argument for NLS with
$k=j=1$. In this case the technicalities are at a minimum and the
essential technique to use oscillations to produce cancellations
is readily apparent. An abbreviated version of the argument in
this section is given on the WEBPAGE of [DZ5]. Finally, in \S4, we
describe our new method as a further development of the steepest
descent method of [DZ1], and use this new method to obtain the
estimates that are needed to prove Theorem 1.10. As in \S3, an
abbreviated version of the estimates in this section is also given
on the WEBPAGE of [DZ5].  Section 4 concludes with a discussion of
the fact, noted at the beginning of this Introduction, that the
error estimate $O(t^{-(1/2+\kappa)})$ in Theorem 1.10 depends only
on the $H^{0,1}$ norm of the initial data $q(x,t=0)=q_0(x)$.
Equation (1.1), however, is not well-posed in $H^{0,1}$. We must
rely instead on the well-known fact (see Theorem 4.39 below) that
equation (1.1) is well-posed in $L^2$, and for $q_0\in H^{0,1}
\subset L^2$ we show that indeed the solution $q=q(x,t)$ obtained
from the RHP (4.1) via (4.3), is precisely the (unique) $L^2$
solution of (1.1) given in Theorem 4.39.

\demo{Remark on Notation (1.11)} If $A=(a_{ij})$ is an $l\times m$
matrix, we define $|A|\equiv
\left(\sum_{i,j}|a_{ij}|^2)\right)^{1/2}=\left(\tr
A^*A\right)^{1/2}$. We say $A(z) =(a_{ij}(z))$ is in $L^p(\Sigma)$
for some contour $\Sigma\subset\cb$ if each of the entries
$a_{ij}(z) \in L^p(\Sigma)$ and we define
$\|A\|_{L^p(\Sigma)}\equiv\|\,|A|\,\|_{L^p(\Sigma)}$.

Throughout the text constants $c>0$ are used generically.
Statements such as $\|f\|\le2c(1+e^c)\le c$, for example, should
not cause any confusion.\enddemo

\subheading{Acknowledgements} The authors would like to thank
Jalal Shatah and Jean Bourgain for useful discussions on the
well-posedness of equation (1.1) in $L^2$. The work of the authors
was supported in part by NSF Grants DMS 0003268 and DMS 0071398
respectively.

\head\S2 Inhomogeneous Riemann-Hilbert Problems\endhead

We begin by summarizing the basic properties of the Cauchy
operator
$$Ch(z) = C_\Gamma h(z) \equiv \int_\Gamma
\frac{h(s)}{s-z} \frac{ds}{2\pi i},\qquad z\in {\Bbb C}\setminus
\Gamma,\tag2.1$$ for an oriented contour $\Gamma$ in the plane. To
fix notation, if we move along the contour in the direction of the
orientation, we say that the $(+)$-side (resp.\ $(-)$-side) lies
to the left (resp.\ right). We have, for example, the following
figure:

\vskip1in
\font\thinlinefont=cmr5

\centerline{\beginpicture \setcoordinatesystem units
<1.00000cm,1.00000cm> \linethickness=1pt \setshadesymbol
({\thinlinefont .}) \setlinear
%
%
\linethickness= 0.500pt \setplotsymbol ({\thinlinefont .})
\putrule from  5.080 19.050 to 12.700 19.050
%
%
\linethickness= 0.500pt \setplotsymbol ({\thinlinefont .})
\putrule from  7.620 21.590 to  7.620 16.510
%
%
\linethickness= 0.500pt \setplotsymbol ({\thinlinefont .})
\putrule from 10.160 21.590 to 10.160 16.510
%
%
\linethickness= 0.500pt \setplotsymbol ({\thinlinefont .})
\putrule from  5.874 19.050 to  6.032 19.050
%
%
\plot  5.779 18.986  6.032 19.050  5.779 19.114 /
%
%
%
\linethickness= 0.500pt \setplotsymbol ({\thinlinefont .})
\putrule from  8.414 19.050 to  8.731 19.050
%
%
\plot  8.477 18.986  8.731 19.050  8.477 19.114 /
%
%
%
\linethickness= 0.500pt \setplotsymbol ({\thinlinefont .})
\putrule from 11.430 19.050 to 11.589 19.050
%
%
\plot 11.335 18.986 11.589 19.050 11.335 19.114 /
%
%
%
\linethickness= 0.500pt \setplotsymbol ({\thinlinefont .})
\putrule from  7.620 21.114 to  7.620 20.796
%
%
\plot  7.557 21.050  7.620 20.796  7.683 21.050 /
%
%
%
\linethickness= 0.500pt \setplotsymbol ({\thinlinefont .})
\putrule from  7.620 17.780 to  7.620 17.462
%
%
\plot  7.557 17.716  7.620 17.462  7.683 17.716 /
%
%
%
\linethickness= 0.500pt \setplotsymbol ({\thinlinefont .})
\putrule from 10.160 21.273 to 10.160 20.955
%
%
\plot 10.097 21.209 10.160 20.955 10.224 21.209 /
%
%
%
\linethickness= 0.500pt \setplotsymbol ({\thinlinefont .})
\putrule from 10.160 17.621 to 10.160 17.780
%
%
\plot 10.224 17.526 10.160 17.780 10.097 17.526 /
%
%
%
\put{$+$} [lB] at  6.191 19.209
%
%
\put{$-$} [lB] at  6.191 18.733
%
%
\put{$+$} [lB] at  9.049 19.209
%
%
\put{$-$} [lB] at  9.049 18.574
%
%
\put{$+$} [lB] at 11.906 19.209
%
%
\put{$-$} [lB] at 11.906 18.733
%
%
\put{$+$} [lB] at  7.779 20.479
%
%
\put{$-$} [lB] at  7.303 20.479
%
%
\put{$+$} [lB] at  7.779 17.145
%
%
\put{$-$} [lB] at  7.303 17.145
%
%
\put{$+$} [lB] at 10.319 20.637
%
%
\put{$-$} [lB] at  9.842 20.637
%
%
\put{$+$} [lB] at  9.842 16.828
%
%
\put{$-$} [lB] at 10.319 16.828 \linethickness=0pt \putrectangle
corners at  5.055 21.615 and 12.725 16.485
\endpicture}

\centerline{Figure 2.2}
\bigskip

\n The properties and estimates that we present below are true for
a very general class of contours, which certainly includes
contours that are finite unions of smooth curves in
$\overline{\Bbb C}$. In particular, the results are true for
contours which are finite unions of straight lines, which is all
that is needed for this paper. We refer the reader to standard
texts on the subject, such as [Dur], for the proofs of the results
that follow.

Suppose $\Gamma$ is given.
\bigskip
\n (2.3)~~Let $h\in L^p(\Gamma, |dz|)$, $1\le p<\infty$. Then
$$C^\pm h(z)\equiv \lim_{\Sb z'\to z\\ z'\in (\pm)
\text{-side of } \Gamma\endSb} (Ch)(z')$$ exists as a
non-tangential limit for a.e.\ $z\in \Gamma$. As usual,
``non-tangential" means that $z' \to z$ in any fixed (truncated)
cone $C$ based at $z$ with $\overline C\setminus\{z\}$ lying
entirely on ${\Bbb C}\setminus\Gamma$.
\bigskip
\n (2.4)~~Let $h\in L^p(\Gamma ,|dz|), 1<p<\infty$. Then
$$\|C^\pm h\|_{L^p(\Gamma)} \le c_p \|h\|_{L^p(\Gamma)
}$$ for some constant $c_p$. Moreover
$$C^\pm h = \pm \frac12 h - \frac12 (Hf)$$
where
$$Hh(z) \equiv \lim_{\vp\downarrow 0} \int\limits_{\Sb\Gamma
\\ |s-z|>\vp\endSb} \frac{h(s)}{z-s} \frac{ds}{i\pi},
\qquad z\in \Gamma,$$ is the Hilbert transform. The above limit
exists a.e.\ on $\Gamma$ and also in $L^p(\Gamma)$, and we have
$$\|Hh\|_{L^p(\Gamma)} \le c'_p \|h\|_{L^p(\Gamma)}$$
for some constant $c'_p$. Clearly
$$C^+-C^- = 1 \quad  \text{and}\quad C^++C^- = -H.$$
(2.5)~~If $\gamma$ is a smooth open segment lying in $\Gamma$,

\vskip1in
\font\thinlinefont=cmr5 \centerline{\beginpicture
\setcoordinatesystem units <1.00000cm,1.00000cm>
\linethickness=1pt
\setshadesymbol ({\thinlinefont .}) \setlinear
%
%
\linethickness= 0.500pt \setplotsymbol ({\thinlinefont .})
\putrule from  2.540 21.590 to  7.620 21.590
%
%
\linethickness= 0.500pt \setplotsymbol ({\thinlinefont .})
\putrule from  6.350 24.130 to  6.350 19.050
%
%
\linethickness= 0.500pt \setplotsymbol ({\thinlinefont .})
\putrule from  9.525 21.590 to 11.430 21.590
%
%
\linethickness= 0.500pt \setplotsymbol ({\thinlinefont .})
\setdashes < 0.1270cm> \plot  7.620 21.590  9.525 21.590 /
%
%
\linethickness= 0.500pt \setplotsymbol ({\thinlinefont .})
\setsolid \plot  8.255 23.019  8.572 21.907 /
%
%
\plot  8.442 22.134  8.572 21.907  8.564 22.169 /
%
%
%
\put{$\gamma$} [lB] at  8.096 23.178 \linethickness=0pt
\putrectangle corners at  2.515 24.155 and 11.455 19.025
\endpicture}

\centerline{Figure 2.6}
\bigskip
\n and $h\in H^1_p(\gamma) = \{f:\ f\in L^p(\gamma, |dz|)$,
$\partial_z f\in L^p(\gamma, |dz|)\}$, $1 < p < \infty$, then
$Ch(z)$ is continuous up to $\gamma$. Thus the limits $C^\pm h(z)$
in (2.3) exists for all $z\in \gamma$ (and without the restriction
of non-tangential convergence). Furthermore, for $h\in
H^1_p(\Gamma)$, $Ch(z)$ is bounded in $\cb\setminus\Gamma$ and
goes to zero uniformly as $z\to\infty$.
\bigskip
\n (2.7)~~It follows, in particular, from (9.3) that if $\Gamma_1$
and $\Gamma_2$ are two contours in ${\Bbb C}$, then for $1 < p <
\infty$
$$\|C^\pm h\|_{L^p(\Gamma_2)} \le c_p\|h\|_{L^p(\Gamma_1)}$$
for $h\in L^p(\Gamma_1)$. For example, if $\Gamma_1= {\Bbb R}$ and
$\Gamma_2 = z_0+{\Bbb R}_+ e^{i\theta}$ for $z_0\in {\Bbb R}$, $0
< \theta <\pi$. Then
$$\|C^\pm h\|_{L^p(z_0+{\Bbb R}_+e^{i\theta})} \le c_p
\|h\|_{L^p({\Bbb R})}$$ for $h\in L^p({\Bbb R})$. Moreover, it is
not hard to see that the constant $c_p$ can be chosen independent
of $z_0\in {\Bbb R}$ and $0 < \theta < \pi$.

Throughout this section we will assume, without further comment,
that $\Gamma$ is a contour for which the above estimates are true.
We now present some of the basic facts about inhomogeneous RHP's.
Many of the results are well known and can be found, implicitly or
explicitly, in standard texts on the subject (see, for example,
[CG]). Other results are new and are tailored specifically to the
analysis of RHP's with external parameters, of the kind that will
arise further on in the paper. For the convenience of the reader
who is interested in applying the theory of inhomogeneous RHP's in
different situations, we present the theory in greater generality
than is needed in the present text.

Consider an oriented contour $\Sigma \subset {\Bbb C}$ as above
with a $k\times k$ {\it jump matrix} $v$: as a standing assumption
throughout the text, we  always assume that $v,v^{-1} \in
L^\infty(\Sigma\to \text{GL}(k,\cb))$. We associate to the pair
$\Sigma, v$ two inhomogeneous RHP's as follows.

For $1<p<\infty$, let $Ch(z) = C_\Sigma h(z) = \int_\Sigma
\frac{h(s)}{s-z} \frac{ds}{2\pi i}$ be the Cauchy operator on
$\Sigma$. Then we say that a pair of $L^p(\Sigma)$-functions
$f_\pm \in \partial C(L^p)$ if there exists a (unique) function
$h\in L^p(\Sigma)$ such that
$$f_\pm(z) = (C^\pm h)(z), \qquad z\in \Sigma. \tag2.8$$ In turn we
will call $f(z) = Ch(z)$, $z\in {\Bbb C}\backslash \Sigma$, the
{\sl extension\/} of $f_\pm = C^\pm h \in \partial C(L^p)$ off
$\Sigma$.

\subhead Inhomogeneous RHP of the first kind IRHP1$
_{L^p}$\endsubhead

Fix $1<p<\infty$. Given $\Sigma, v$ and a function $f$, we say
that $m_\pm \in f + \partial C(L^p)$ solves an IRHP1$_{L^p}$ if
$$m_+(z) = m_-(z) v(z),\qquad z\in \Sigma.$$

\subhead Inhomogeneous RHP of the second kind
IRHP2$_{L^p}$\endsubhead

Fix $1<p<\infty$. Given $\Sigma, v$ and a function $F\in
L^p(\Sigma)$, we say that
 $M_\pm \in \partial C(L^p)$ solves an IRHP2$_{L^p}$ if
$$M_+(z) = M_-(z) v(z) + F(z),\qquad z\in \Sigma.$$

Recall that $m$ solves the normalized RHP $(\Sigma, v)$ if, at
least
 formally (cf. {1.4) above),

 $\bullet \ m(z)$ is analytic in ${\Bbb C}\backslash \Sigma$

$\bullet \ m_+(z) = m_-(z) v(z), z\in \Sigma$

$\bullet \ m(z) \to I$ as $z\to \infty$.
\medskip

\n More precisely, we make the following definition.

\definition{Definition 2.9}
Fix $1<p<\infty$. We say that $m_\pm$ solves the normalized RHP
$(\Sigma,v)_{L^p}$ if
 $m_\pm$ solves the IRHP1$_{L^p}$ with $f\equiv I$.

In the above definition, if $m_\pm-I= C^\pm h,$ then clearly the
extension $m=I+ Ch$ of $m_\pm$ off $\Sigma$ solves the above RHP
in the formal sense. If $p=2$, which is mostly the case of
interest, we will drop the subscript and simply write  IRHP1,
IRHP2 and $(\Sigma,v)$. Let $v = (v^-)^{-1} v^+ = (I-w^-)^{-1}
(I+w^+)$ be a factorization of $v$ with
 $v^\pm,(v^{\pm})^{-1}\in L^\infty$, and let $C_w, w = (w^-,w^+)$, denote the basic,
associated singular integral operator
$$C_w h = C^+(hw^-) + C^-(hw^+)
\tag2.10$$ acting on $L^p(\Sigma)$-matrix-valued functions $h$. As
$w^\pm\in L^\infty$, $C_w$ is clearly bounded from $L^p\to L^p$
for all $1<p<\infty$. The utility of IRHP1$_{L^p}$ and
IRHP2$_{L^p}$ lies in the following circle of ideas. As we will
shortly see, the solution $m_\pm$ of an IRHP1$_{L^p}$ is directly
related to the inverse of $1-C_w$. However, for general $f$ (as
opposed to the case $f\equiv I$ in the normalized RHP above),
$m_\pm$ does not have an analytic continuation to ${\Bbb
C}\backslash \Sigma$. On the other hand, the solution $M_\pm$ of
an IRHP2$_{L^p}$ corresponds to the boundary values of an analytic
function. The key point, as we shall see, is that under
appropriate assumptions on $v$ and $f$, IRHP1$_{L^p}$ and
IRHP2$_{L^p}$ are equivalent, and hence steepest-descent type
methods applied directly to $M_\pm$ can be used to control the
resolvent $(1-C_w)^{-1}$.

Now suppose $f$ and $v$ are such that $f(v-I)\in L^p(\Sigma)$ for
some $1<p<\infty$. Let $M_\pm$ solve the IRHP2$_{L^p}$ with $F =
f(v-I)$. Then simple algebra shows that $m_\pm = M_\pm + f$ solves
the IRHP1$_{L^p}$, $m_+ = m_- v$, $m_\pm - f \in \partial C(L^p)$.
If we reverse the argument, we see that if $m_\pm$ solves the
IRHP1$_{L^p}$ with any given $f$, then $M_\pm = m_\pm -f$ solves
IRHP2$_{L^p}$ with $F$ of the form $f(v-I)$ which, however, is not
the general case as $v-I$ need not be invertible. To prove the
equivalence of IRHP1$_{L^p}$ and IRHP2$_{L^p}$, we must proceed in
a different way. Given $F\in L^p(\Sigma)$, let $m_\pm$ solve
IRHP1$_{L^p}$ with $f= C^-F$. By assumption, $m_\pm - f = C^\pm h$
for some $h\in L^p(\Sigma)$, and hence
$$\align m_-  &= f + C^-h = C^-(F+h),\\ m_+ + F &= f + C^+ h + (C^+
- C^-)F = C^+(F+h). \endalign$$ It follows that $M_+ = m_+ +F$,
$M_- = m_-$ solve IRHP2$_{L^p}$ with the given $F$. We have proved
the following result.
\enddefinition

\proclaim{Proposition 2.11 \rom{(Equivalence of IRHP1$_{L^p}$ and
IRHP2$_{L^p}$)}} Suppose $f$ and $v$ are such that $f(v-I) \in
L^p(\Sigma)$ for some $1<p<\infty$. Then
$$m_\pm = M_\pm +f\tag2.12$$
solves IRHP1$_{L^p}$ with the given $f$, if $M_\pm$ solve
IRHP2$_{L^p}$ with $F = f(v-I)$. Conversely, if $F\in
L^p(\Sigma)$, then
$$
M_+ = m_+ +F, \ \ M_- = m_- \tag2.13$$ solve IRHP2$_{L^p}$ with
the given $F$, if $m_\pm$ solve IRHP1$_{L^p}$ with $f =
C^-F$.\endproclaim

We will be interested principally in the situation where $f \in
 L^\infty(\Sigma)$ and $v-I \in L^p(\Sigma)$ for some $1<p<\infty$.

To establish the connection with the inverse of $1-C_w$, let
$m_\pm$
 solve IRHP1$_{L^p}$ with\break $f\in L^p(\Sigma)$. Thus $m_\pm = f + C^\pm h$
for some $h\in L^p(\Sigma)$. Also $m_+ = m_-v \equiv
m_-(v^-)^{-1}v^+$. Set $\mu = m_+(v^+)^{-1} = m_-(v^-)^{-1} \in
L^p(\Sigma)$, and define $H(z) = (C(\mu (w^+ + w^-))(z)$, $z\in
{\cb} \setminus\Sigma$. Then a simple calculation shows that
$H_\pm = (C_w-1) \mu + m_\pm$, or $m_\pm - f - H_\pm = (1-C_w)\mu
- f$. But $m_\pm - f - H_\pm \in \partial   C(L^p)$, and it
follows that $(1-C_w)\mu = f$. Conversely, if $\mu\in L^p(\Sigma)$
solves $(1-C_w)\mu = f$, then $H = C(\mu(w^+ + w^-))$ satisfies
$H_\pm = -f + \mu v^\pm$. Thus setting $m_\pm = \mu v^\pm$, we see
that $m_+ = m_-v$ and $m_\pm - f \in \partial   C(L^p)$. In
particular, $\mu\in L^p$ solves $(1-C_w) \mu = 0$ if and only if
$m_\pm = \mu v_\pm$ solves the homogeneous RHP $$m_+ = m_-v,\quad
m_\pm \in
\partial   C(L^p).$$ We summarize the above calculations
as follows.

\proclaim{Proposition 2.14} Let $1<p<\infty$. Then
$$\gather
1-C_w \text{ is a bijection in } L^p(\Sigma)\\
\Leftrightarrow\\
\text{IRHP1$_{L^p}$ has a unique solution for all } f \in L^p(\Sigma)\\
\Leftrightarrow\\
\text{IRHP2$_{L^p}$ has a unique solution for all } F\in
L^p(\Sigma).
\endgather$$
Moreover, if $(1-C_w)^{-1}$ exists, then for $f\in L^p(\Sigma)$,
$$\align
(1-C_w)^{-1} f &= m_+(v^+)^{-1} = m_-(v^-)^{-1}\tag2.15\\
&= (M_+ + f) (v^+)^{-1} = (M_- + f) (v^-)^{-1}
\endalign$$
where $m_\pm$ solve IRHP1$_{L^p}$ with the given $f$, and $M_\pm$
solve IRHP2$_{L^p}$ with $F = f(v-I)$. Conversely, if $M_\pm$
solves IRHP2$_{L^p}$ with $F\in L^p(\Sigma)$, then $$ M_+ =
((1-C_w)^{-1} C^-F) v^++F \text{ \ \ and \ \ } M_- = ((1-C_w)^{-1}
C^-F)v^-.\tag2.16$$ Finally, if $f \in L^\infty(\Sigma)$ and
$v^\pm -I \in L^p(\Sigma)$, then (2.15) remains valid provided we
interpret
$$(1-C_w)^{-1} f \equiv f + (1-C_w)^{-1} (C_w f). \tag2.17$$
This is true in particular for the normalized RHP $(\Sigma, v)_p$
where $f\equiv I$.
\endproclaim
The above Proposition implies, in particular, that if $\mu\in I+
L^p(\Sigma)$ solves
$$(1-C_w)\mu=I\tag2.18$$ as in (2.17), then $m_\pm=\mu v^\pm$ solves
the normalized  RHP (2.9).

The fact that IRHP1$_{L^p}$ and IRHP2$_{L^p}$ depend on $v$, and
not on the particular factorization $v = (v^-)^{-1} v^+$, has the
following immediate consequence.

\proclaim{Corollary to Proposition 2.14} Let as $1<p<\infty$. The
operator $1-C_w$ is bijective in $L^p(\Sigma)$ for all
factorizations $v = (v^-)^{-1} v^+ = (I-w^-)^{-1} (I+w^+)$, if and
only if $1-C_{w'}$ is bijective for at least one factorization $v
= (v'{}^-)^{-1} v'{}^+ = (I-w'{}^-)^{-1} (I-w'{}^+)$. Moreover,
for $f\in L^p(\Sigma)$ $$(1-C_w)^{-1} f = [(1-C_{w'})^{-1} f]b$$
where $b = v'{}^+(v^+)^{-1}  = v'{}^-(v^-)^{-1}$. \endproclaim

We now describe a useful and simple fact about the operator $C_w$.
Given $\Sigma$ and $v$, suppose we reverse the orientation of
$\Sigma$ on some subset $\gamma\subset \Sigma$. For example,
consider the following figure.
\bigskip\bigskip
\font\thinlinefont=cmr5 \centerline{\beginpicture
\setcoordinatesystem units <0.7cm,0.7cm>
\linethickness=1pt
\setshadesymbol ({\thinlinefont .}) \setlinear
%
%
\linethickness= 0.500pt \setplotsymbol ({\thinlinefont .})
\putrule from  2.540 21.590 to  7.620 21.590
%
%
\linethickness= 0.500pt \setplotsymbol ({\thinlinefont .})
\putrule from 12.700 21.590 to 17.780 21.590
%
%
\linethickness= 0.500pt \setplotsymbol ({\thinlinefont .})
\putrule from  3.810 24.130 to  3.810 19.050
%
%
\linethickness= 0.500pt \setplotsymbol ({\thinlinefont .})
\putrule from  6.350 24.130 to  6.350 19.050
%
%
\linethickness= 0.500pt \setplotsymbol ({\thinlinefont .})
\putrule from 13.970 24.130 to 13.970 19.050
%
%
\linethickness= 0.500pt \setplotsymbol ({\thinlinefont .})
\putrule from 16.351 24.130 to 16.351 19.050
%
%
\linethickness= 0.500pt \setplotsymbol ({\thinlinefont .})
\putrule from  9.207 21.590 to 10.795 21.590
%
%
\plot 10.541 21.526 10.795 21.590 10.541 21.654 /
%
%
%
\linethickness= 0.500pt \setplotsymbol ({\thinlinefont .})
\putrule from  2.699 21.590 to  3.016 21.590
%
%
\plot  2.762 21.526  3.016 21.590  2.762 21.654 /
%
%
%
\linethickness= 0.500pt \setplotsymbol ({\thinlinefont .})
\putrule from  4.763 21.590 to  4.921 21.590
%
%
\plot  4.667 21.526  4.921 21.590  4.667 21.654 /
%
%
%
\linethickness= 0.500pt \setplotsymbol ({\thinlinefont .})
\putrule from  6.985 21.590 to  7.144 21.590
%
%
\plot  6.890 21.526  7.144 21.590  6.890 21.654 /
%
%
%
\linethickness= 0.500pt \setplotsymbol ({\thinlinefont .})
\putrule from 13.018 21.590 to 13.176 21.590
%
%
\plot 12.922 21.526 13.176 21.590 12.922 21.654 /
%
%
%
\linethickness= 0.500pt \setplotsymbol ({\thinlinefont .})
\putrule from 14.764 21.590 to 14.922 21.590
%
%
\plot 14.669 21.526 14.922 21.590 14.669 21.654 /
%
%
%
\linethickness= 0.500pt \setplotsymbol ({\thinlinefont .}) \plot
17.145 21.590 17.145 21.590 /
%
%
\linethickness= 0.500pt \setplotsymbol ({\thinlinefont .})
\putrule from 16.986 21.590 to 17.145 21.590
%
%
\plot 16.891 21.526 17.145 21.590 16.891 21.654 /
%
%
%
\linethickness= 0.500pt \setplotsymbol ({\thinlinefont .})
\putrule from  3.810 20.161 to  3.810 20.320
%
%
\plot  3.873 20.066  3.810 20.320  3.747 20.066 /
%
%
%
\linethickness= 0.500pt \setplotsymbol ({\thinlinefont .})
\putrule from  6.350 19.844 to  6.350 20.003
%
%
\plot  6.413 19.748  6.350 20.003  6.287 19.748 /
%
%
%
\linethickness= 0.500pt \setplotsymbol ({\thinlinefont .})
\putrule from 13.970 20.479 to 13.970 20.320
%
%
\plot 13.906 20.574 13.970 20.320 14.034 20.574 /
%
%
%
\linethickness= 0.500pt \setplotsymbol ({\thinlinefont .}) \plot
16.351 20.637 16.351 20.637 /
%
%
\linethickness= 0.500pt \setplotsymbol ({\thinlinefont .})
\putrule from 16.351 20.796 to 16.351 20.637
%
%
\plot 16.288 20.892 16.351 20.637 16.415 20.892 /
%
%
%
\linethickness= 0.500pt \setplotsymbol ({\thinlinefont .})
\putrule from  3.810 22.701 to  3.810 22.860
%
%
\plot  3.873 22.606  3.810 22.860  3.747 22.606 /
%
%
%
\linethickness= 0.500pt \setplotsymbol ({\thinlinefont .})
\putrule from  6.350 22.384 to  6.350 22.543
%
%
\plot  6.413 22.289  6.350 22.543  6.287 22.289 /
%
%
%
\linethickness= 0.500pt \setplotsymbol ({\thinlinefont .})
\putrule from 13.970 22.860 to 13.970 23.019
%
%
\plot 14.034 22.765 13.970 23.019 13.906 22.765 /
%
%
%
\linethickness= 0.500pt \setplotsymbol ({\thinlinefont .})
\putrule from 16.351 22.543 to 16.351 22.701
%
%
\plot 16.415 22.447 16.351 22.701 16.288 22.447 /
%
%
%
\linethickness= 0.500pt \setplotsymbol ({\thinlinefont .}) \plot
14.764 19.526 14.287 20.161 /
%
%
\plot 14.491 19.996 14.287 20.161 14.389 19.920 /
%
%
%
\linethickness= 0.500pt \setplotsymbol ({\thinlinefont .}) \plot
15.399 19.526 16.034 20.161 /
%
%
\plot 15.899 19.937 16.034 20.161 15.809 20.027 /
%
%
%
\put{$\gamma$} [lB] at 14.922 19.209 \linethickness=0pt
\putrectangle corners at  2.515 24.155 and 17.805 19.025
\endpicture}
\centerline{Figure \ \  2.19}
\bigskip

\n If we denote the new contour by $\widehat\Sigma$ and set $\hat
v \equiv v$ on $\Sigma \setminus \gamma$, $\hat v \equiv v^{-1}$
on $\gamma$, then we clearly have equivalent RHP's on $\Sigma,v$
and on $\widehat\Sigma, \hat v$. Furthermore, a factorization $v =
(v^-)^{-1}v^+$, leads to the natural factorization $\hat v =
(v^-)^{-1} v^+$ on $\Sigma\setminus \gamma$ and $\hat v =
(v^+)^{-1}v^-$ on $\gamma$. Thus, where $v^\pm = I \pm w^\pm$,
$\hat v^\pm = I \pm \widehat w^\pm$, we have $\widehat w^\pm =
w^\pm$ on $\Sigma\setminus \gamma$, $\widehat w^\pm = -w^\mp$ on
$\gamma$.

\proclaim{Proposition 2.20} Let $\widehat\Sigma, \hat v$ be
defined as above. Then
$$C_{\widehat w} = C_w\tag2.21$$
in $L^p(\Sigma) = L^p(\widehat\Sigma)$ for any $1<p<\infty$.
\endproclaim

\demo{Proof} Let $h\in L^p(\Sigma) = L^p(\widehat \Sigma)$. Then
it follows directly from the above definitions that $C_\Sigma
(h(w^+ +w^-)) = C_{\widehat\Sigma} (h(\widehat w^+ + \widehat
w^-))$, where $C_\Sigma, C_{\widehat\Sigma}$ denote the Cauchy
operators on the (oriented) contours $\Sigma, \widehat\Sigma$
respectively. Now
$$C^\pm_\Sigma (h(w^++w^-)) = C_wh \pm hw^\pm, \quad C^\pm_{\widehat\Sigma} (h(\widehat w^+ + \widehat w^-)) = C_{\widehat w} h \pm h\widehat w^\pm.$$
On $\Sigma\setminus \gamma, C^\pm_\Sigma (h(w^+ + w^-)) =
C^\pm_{\widehat\Sigma} (h(\widehat w^+ + \widehat w^-))$ and
$w^\pm = \widehat w^\pm$, and it follows that $C_wh = C_{\widehat
w}h$. But on $\gamma$, $C^\pm_\Sigma(h(w^++w^-))$
$=C^\mp_{\widehat\Sigma} (h(\widehat w^+ + \widehat w^-))$ and
$w^\pm = -\widehat w^\mp$, and it follows again that $C_wh =
C_{\widehat w}h$. The proof is done. $\qquad\square$\enddemo

Finally we consider uniqueness for the solution of the normalized
RHP $(\Sigma,v)_p$ as given in Definition~2.9. Observe first that
if $F(z) = (Cf)(z)$ for $f\in L^p(\Sigma)$ and $G(z) = (Cg)(z)$
for $g\in L^q(\Sigma)$, $\frac1r=\frac1p + \frac1q\le1$,
$1<p,q<\infty$, then a simple computation shows that
$$(FG)(z) =(Ch)(z)\tag2.22$$
 where
$$h(s) = - \frac12 (g(s) (Hf)(s) + f(s)
(Hg)(s))\tag2.23$$ and $Hf(s) = \lim\limits_{\vp \downarrow 0}
\int\limits_{|s'-s|>\vp} \frac{f(s')}{s-s'} \frac{ds'}{i\pi}$,
$Hg(s) = \lim\limits_{\vp\downarrow 0} \int\limits_{|s'-s|>\vp}
\frac{g(s')}{s-s'} \frac{ds'}{i\pi}$ denote the Hilbert transforms
of $f$ and $g$ respectively (see (2.4)). As $h$ clearly lies in
$L^r(\Sigma)$, it follows that
$$F_+(z) G_+(z) - F_-(z) G_-(z) = h(z)\tag2.24$$
for a.e.\ $z\in \Sigma$. It is now easy to prove uniqueness for
the normalized RHP in the following form.

\proclaim{Theorem 2.25} Fix $1<p<\infty$. Suppose $m(z)$ solves
the normalized RHP $(\Sigma,v)_{L^p}$
 in the sense of Definition~2.9, and suppose in addition that $m(z)^{-1}$
exists for all $z\in {\Bbb C}\setminus \Sigma$ and that
$(m^{-1})_\pm - I \in \partial   C(L^q)$, $1<q<\infty$, $\frac1r
\equiv \frac1p + \frac1q \le 1$. Then $m$ is unique.
\endproclaim

\demo{Proof} Suppose $\widehat m$ is a second solution of the
normalized RHP. Then arguing as above, $\widehat m m^{-1} - I =
(\widehat m-I)(m^{-1}-I) + (\widehat m-I) + (m^{-1}-I) = Ch$ where
$h\in L^r(\Sigma) + L^p(\Sigma) + L^q(\Sigma)$, and $\widehat m_+
m^{-1}_+ - \widehat m _- m^{-1}_- = h$. But $\widehat m_+m^{-1}_+
= (\widehat m_-v) (m_-v)^{-1}  = \widehat m_- m^{-1}_-$ a.e.\ on
$\Sigma$, and hence $h=0$. Thus $\widehat mm^{-1} -I = 0$, or
$\widehat m=m$. $\qquad\square$\enddemo

\remark{Remark} Observe that in the above proof we did not need to
assume that also $(\widehat m^{-1})_\pm - I \in \partial C(L^q)$.

In the special case in which $k=2$, $p=2$ and $\det v(z) = 1$
a.e.\ on $\Sigma$, there is the following stronger result.
\endremark

\proclaim{Theorem 2.26} If $k=2, p=2$ and $\det v(z) =1$ a.e.\ on
$\Sigma$, then the solution $m$ of
 the normalized RHP $(\Sigma, v) = (\Sigma,v)_{L^2}$ is unique.
\endproclaim

\demo{Proof} As $k=2$ and $p=2$, (2.22) and (2.23) imply that
$\det m(z) - 1=(Ch)(z)$ where $h\in L^1(\Sigma)+L^2(\Sigma)$, and
$\det m_+(z) - \det m_-(z) = h(z)$ a.e.\ on $\Sigma$. But $\det
m_+(z) = \det m_-(z) \det v(z)$ $=\det m_-(z)$, and hence $h=0$
and so $\det m(z)=1$. But if $m = (m_{ij})_{1\le
i,j\le 2}$, then $m^{-1} = \left(\matrix m_{22}&-m_{12}\\
-m_{21}&m_{11}\endmatrix\right)$ and clearly $(m^{-1})_\pm - I \in
\partial   C(L^2)$. The result then follows from
Theorem~2.25.$\qquad\square$ \enddemo

The above result is of considerable general interest in practice,
and applies, in particular, to the RHP (1.4) for the NLS equation.
But in the case $\Sigma=\rb$ with the special jump matrix $v_x$ in
(1.4) uniqueness and existence for the normalized RHP follows
immediately from the following observation. Write

$$\pmtwo1{-re^{izx}}01^{-1}\pmtwo10{-\bar re^{-izx}}1=
(I-w_x^-)^{-1}(I+w_x^+)\tag2.27$$so that$$(w_x^-,w_x^+)=\left(
\pmtwo0{re^{izx}}00,\pmtwo00{-\bar re^{-izx}}0\right).\tag2.28$$
Then for $C_{w_x}$ in (2.10). and $h=(h_{ij})_{1\le i,j\le2}$,
$$C_{w_x}h=\pmtwo{C^-(-h_{12}\bar re^{-i\loz x})}{C^+(h_{11} re^{i\loz x})}
{C^-(-h_{22}\bar re^{-i\loz x})}{C^+(h_{21} re^{i\loz x})}.$$ But,
under the Fourier transform, $C^+ $ (resp.$-C^-$) is just
multiplication by the characteristic function of $(0,\infty)$
(resp. $(-\infty,0)$) and hence $\|C^\pm\|_{L^2(\rb)}=1$. Thus
(recall Remark (1.10)), $\norm2{C_{w_x}}\le\norm\infty r\norm2h$
and hence $$\norm2{C_{w_x}}\le\norm\infty r.\tag2.29$$

As noted in the Introduction, we always have $\norm\infty r<1$,
and it follows that $(1-C_{w_x})^{-1}$ exists and is uniformly
bounded on $L^2$ for all
$x\in\rb$,$$\norm2{(1-C_{w_x})^{-1}}\le\frac1{1-\norm\infty
r}\tag2.30$$ By Section 3, we always have $r\in L^2(\rb)$. It
follows then by Proposition 2.14 that we have proved the following
result.

\proclaim{Proposition 2.31} The normalized RHP $(\Sigma=\rb,v_x)$
has a unique solution $m_\pm\in I+\dc2$ for each $x\in\rb$.
Moreover $$m_\pm=I+C^\pm(\mu(w_x^++w_x^-)),\tag2.32$$where $\mu\in
I+L^2(\rb)$ is the unique solution of
$$(1-\cxs)\mu=I.\tag2.33$$\endproclaim

In the remainder of this section, we use the relationships between
various inhomogeneous RHP's to establish a uniform bound (see
Proposition 2.49 below) of type (2.30) for an associated model RHP
which plays a key role in proving Theorem 1.10. We refer the
reader to \S4 for an explanation of the origin and relevance of
the model.

By (1.7)(1.8) we need to consider the normalized RHP (1.4) with a
time-dependent reflection coefficient $r(t)=e^{-iz^2t}r$ and hence
a space and time dependent jump matrix
$$v_\theta=\pmatrix1-|r(z)|^2&r(z)e^{i\theta}\\-\overline{r(z)}e^{-i\theta}&1\endpmatrix,\tag2.34$$
where $$\theta=xz-tz^2.\tag2.35$$ Throughout the paper, if $A$ is
a $2\times2$ matrix,
 then $A_\theta$ denotes $e^{i\theta \ads}A\equiv e^{i\theta\sigma}Ae^{-i\theta\sigma}$, where $
\sigma=\diag(1/2,-1/2)$ as before. In this notation,
$v_\theta=e^{i\theta \ads}v$, where
$v=\pmatrix1-|r(z)|^2&r(z)\\-\overline{r(z)}&1\endpmatrix.$

Let $z_0=x/(2t)$ denote the stationary phase point for
$e^{i\theta}$. For a function $f$ on $\rb$ we introduce the
notation $$[f](z)=\frac {f(z_0)}{(1+i(z-z_0))^2},\ \ \
z\in\rb.\tag2.36$$ Clearly $[f](z_0)=f(z_0).$

Let $\delta_\pm$ be the solution of the scalar, normalized RHP
$({\Bbb R}_- + z_0, 1- |r|^2)$,
$$\left\{\aligned
&\delta_+ = \delta_- (1-|r|^2),\qquad z\in {\Bbb R}_-
+z_0,\\
&\delta_\pm - 1 \in \partial   C (L^2),
\endaligned\right.\tag2.37$$
where the contour ${\Bbb R}_- + z_0$ is oriented from $-\infty$ to
$z_0$. The properties of $\delta$ can be read off from the
following elementary proposition, whose proof is left to the
reader.

\proclaim{Proposition 2.38} Suppose $r\in L^\infty({\Bbb R}) \cap
L^2({\Bbb R})$ and $\|r\|_{L^\infty} \le \rho <1$. Then the
solution $\delta_\pm$ of the scalar normalized RHP (2.37) exists
and is unique and is given by the formula
$$\delta_\pm(z) = e^{C^\pm_{{\Bbb R}_-+z_0} \log
1-|r|^2} = e ^{\frac1{2\pi i} \int^{z_0}_{-\infty}
\frac{\log(1-|r(s)|^2)}{s-z_\pm} ds},\qquad z\in {\Bbb
R}.\tag2.39$$ The extension $\delta$ of $\delta_\pm$ off ${\Bbb
R}_- +z_0$ is given by
$$\delta(z) = e^{C_{{\Bbb R}_- +z_0} \log(1-|r|^2)}
= e^{\frac1{2\pi i} \int^{z_0}_{-\infty} \frac{\log
(1-|r(s)|^2)}{s-z} ds},\qquad z\in {\Bbb C} \setminus ({\Bbb R}_-
+ z_0),\tag2.40$$ and satisfies for $z\in {\Bbb C}\setminus ({\Bbb
R}_- +z_0)$,
$$\delta(z) \overline{\delta(\bar z)} = 1,\tag2.41$$
$$(1-\rho)^{\frac12} \le (1-\rho^2)^{\frac12} \le
|\delta(z)|, |\delta^{-1}(z)| \le (1-\rho^2)^{-\frac12} \le
(1-\rho)^{-\frac12},\tag2.42$$ and
$$|\delta^{\pm 1}(z)|\le 1 \quad \text{for}\quad
\pm \text{Im } z>0.\tag2.43$$

For real $z$,
$$|\delta_+(z) \delta_-(z)| = 1 \quad \text{and, in
particular,}\quad |\delta(z)|  = 1\quad \text{for} \quad
z>z_0,\tag2.44$$
$$|\delta_+(z)| = |\delta^{-1}_-(z)| = (1 - |r(z)|^2
)^{\frac12},\qquad z<z_0,\tag2.45$$ and
$$\Delta\equiv \delta_+\delta_- = e^{\frac1{i\pi}
\text{ P.V. } \int^{z_0}_{-\infty}
\frac{\log(1-|r(s)|^2)}{s-z}ds}, \quad \text{where P.V.\ denotes
the principal value.} \tag2.46$$ Also $|\Delta| =
|\delta_+\delta_-| = 1$.
$$\|\delta_\pm - 1\|_{L^2(dz)} \le \frac{c\|r\|_{L^2}}
{1-\rho}.\tag2.47$$\endproclaim

For a reflection coefficient $r(z)$, $\norm\infty r<1$, set
$$\hat w=(\hat w^-,\hat w^+)=\cases\left(\pmtwo0{[r]\delta^2}00,\pmtwo00{-[\bar
r]\delta^{-2}}0\right),\ \ \
&z>z_0,\\\left(\pmtwo00{-\left[\frac{\bar
r}{1-|r|^2}\right]\delta_-^{-2}}0,\pmtwo0{\left[\frac
r{1-|r|^2}\right]\delta_+^2}00\right),\ \ \
&z<z_0,\endcases\tag2.48$$and as above, set $\hat
w_\theta=\epts\hat w=(\epts\hat w^-,\epts\hat w^+).$ We consider
the model normalized RHP $(\Sigma=\rb,\hat v_\theta)$, where $\hat
v_\theta=(I-\hat w_\theta^-)^{-1})(I+\hat w_\theta^+)$. Our goal
is to prove the following analog of the resolvent bound in (2.30).

\proclaim{Proposition 2.49} For $t$ sufficiently large, say $t\ge
t_0$, $(1-\cht)^{-1}$ exists in $L^2(\rb)$ and for some constant
c,$$\norm2{(1-\cht)^{-1}}\le c \tag2.50$$for all $x\in\rb$ and for
all $t\ge t_0$.
\endproclaim
In order to understand the origin of this bound, consider $z$
close to $z_0$, so that $[r](z)\sim r(z)$ etc., a direct
calculation using the jump relation (2.37) shows that $\hat
v_\theta(z)\sim\delta_-^{\sigma_3}v_\theta\delta_+^{-\sigma_3},$
where $v_\theta$ is given in (2.34) and
$\sigma_3=2\sigma=\diag(1,-1)$ is the third Pauli matrix. Now by
(2.30), $(1-\ct)^{-1}$ exists and
$\norm2{(1-\ct)^{-1}}\le\frac1{1-\norm\infty
{r(t)}}=\frac1{1-\norm\infty r}$ for all $x,t$ and $r(t)=e^{-iz^2
t}r$. It follows then from Proposition 2.14 that the IRHP2,
$$M_+=M_- v_\theta+F,\ \ \ M_\pm\in\dc2$$ is solvable and $\norm2{{M_\pm}}\le\frac
c{1-\norm\infty r}\norm2 F$ for all $F\in L^2$, uniformly for all
$x,t$. But then $M_\pm^\delta\equiv M_\pm\delta_\pm^{-\sigma_3}$
solves
$$M_+^\delta=M_-^\delta(\delta_-^{\sigma_3}v_\theta\delta_+^{-\sigma_3})+F^\delta$$
with $F^\delta=F\delta_+^{-\sigma_3}$. By the analyticity and
boundedness properties of $\delta(z)$ in Proposition 2.38, it is
easy to see (cf. the approximation argument $h^\ep\to h$ preceding
(2.56) below) that $M_\pm^\delta\in\dc2$ and $F^\delta\in L^2$. As
$\delta_\pm,\delta_\pm^{-1}$ are bounded, what the above
calculation in fact shows is that the IRHP2 for
$\delta_-^{\sigma_3}v_\theta\delta_+^{-\sigma_3}$ is uniquely
solvable for all $F^\delta\in L^2$ and $\|M_\pm^\delta\|\le\frac
c{1-\norm\infty r}\norm2{F^\delta}$ uniformly for all $x,t$. Again
by Proposition 2.14 we infer finally the analog of the uniform
bound (2.50) for
$\delta_-^{\sigma_3}v_\theta\delta_+^{-\sigma_3}$. The thrust of
the argument that follows is that as $t\to\infty$ the RHP
$(\rb,\hat v_\theta)$ indeed ``localizes" near $z=z_0$ and
$\delta_-^{\sigma_3}v_\theta\delta_+^{-\sigma_3}$ is a good enough
approximation to $\hat v_\theta(z)$ to infer the desired bound on
$(1-\cht)^{-1}$.

Consider the extended problem on
$\Sigma^e=\rb\cup(z_0+e^{i\pi/4}\rb)\cup(z_0+e^{-i\pi/4}\rb)$
oriented as in Figure 2.51

\bigskip\bigskip
\font\thinlinefont=cmr5

\centerline{\beginpicture \setcoordinatesystem units <.7cm,.7cm>
\setshadesymbol ({\thinlinefont .}) \setlinear
%
%
\linethickness= 0.500pt \setplotsymbol ({\thinlinefont .}) \plot
2.540 24.130 12.700 19.050 /
%
%
\linethickness= 0.500pt \setplotsymbol ({\thinlinefont .}) \plot
2.540 19.050 12.700 24.130 /
%
%
\linethickness= 0.500pt \setplotsymbol ({\thinlinefont .})
\putrule from  2.540 21.590 to 12.700 21.590
%
%
\linethickness= 0.500pt \setplotsymbol ({\thinlinefont .})
\putrule from  3.810 21.590 to  4.128 21.590
%
%
\plot  3.873 21.526  4.128 21.590  3.873 21.654 /
%
%
%
\linethickness= 0.500pt \setplotsymbol ({\thinlinefont .})
\putrule from 10.160 21.590 to 10.478 21.590
%
%
\plot 10.224 21.526 10.478 21.590 10.224 21.654 /
%
%
%
\linethickness= 0.500pt \setplotsymbol ({\thinlinefont .}) \plot
3.810 19.685  4.128 19.844 /
%
%
\plot  3.929 19.673  4.128 19.844  3.872 19.787 /
%
%
%
\linethickness= 0.500pt \setplotsymbol ({\thinlinefont .}) \plot
3.493 23.654  3.810 23.495 /
%
%
\plot  3.554 23.552  3.810 23.495  3.611 23.665 /
%
%
%
\linethickness= 0.500pt \setplotsymbol ({\thinlinefont .}) \plot
10.160 22.860 10.478 23.019 /
%
%
\plot 10.279 22.848 10.478 23.019 10.222 22.962 /
%
%
%
\linethickness= 0.500pt \setplotsymbol ({\thinlinefont .}) \plot
9.842 20.479 10.160 20.320 /
%
%
\plot  9.904 20.377 10.160 20.320  9.961 20.490 /
%
%
%
\put{$z_0$} [lB] at  7.461 20.955
%
%
%
%
%
%
\put{$\pi/4$} [lB] at  5.5 21.8
%
%
\put{$\pi/4$} [lB] at 8.75 21.8
%
%
\put{$\pi/4$} [lB] at  5.5 21.1
%
%
\put{$\pi/4$} [lB] at 8.75  21.1 \linethickness=0pt \putrectangle
corners at  2.515 24.155 and 12.725 19.025
\endpicture}

\centerline{Figure 2.51 $\Sigma^e$}
\bigskip

\noindent with $v^e\equiv\hat v_\theta$ on $\rb$ and $v^e\equiv I$
on $\Sigma^e\setminus\rb$. The opening angle $\pi/4$ is chosen
only for convenience;  any angle between 0 and $\pi/2$ would do.

We now show that the bound (2.50) on $(1-\cht)^{-1}$ is equivalent
to the bound  $$\Norm2{\Sigma^e}{(1-\cht)^{-1}}\le c \tag2.52$$
for all $x\in\rb$ and $t\ge t_0$ for some $t_0$. Here $w^e=(\hat
w_\theta^-,\hat w_\theta^+)$ on $\rb$ and $w^e=0$ on
$\Sigma^e\setminus\rb$. Suppose $(1-\cht)^{-1}$ exists in
$L^2(\rb)$ and consider the equation $(1-C_{w^e})\mu^e=F^e$ for
$F^e\in L^2(\Sigma^e)$. Then as $w^{e\pm}=0$ on $
\Sigma^e\setminus\rb$, we obtain a decoupled equation for
$\mu\restriction\rb$ alone,

$$(1-\cht)(\mu^e\restriction\rb)=F^e\restriction\rb$$
on $\rb$, and hence
$$(\mu^e\restriction\rb)=(1-\cht)^{-1}(F^e\restriction\rb),$$
and on $\Sigma^e\setminus\rb)$,
$$\mu^e=F^e+C((\mu^e\restriction\rb)(\hat w_\theta^++\hat
w_\theta^-))$$ and therefore, using the properties of the Cauchy
operator and the fact that $\norm\infty r<1$,
$$\Norm2{\Sigma^e}{\mu^e}\le
c\Norm2\rb{(1-\cht)^{-1}}\Norm2{\Sigma^e}{F^e}.$$Conversely,
suppose $(1-C_{w^e})^{-1}$ exists in $L^2(\Sigma^e)$. Then if
$\hat F\in L^2(\rb)$, we set $F^e=\hat F$ on $\rb$ and zero on
$\Sigma^e\setminus\rb$. Let $\mu^e$ solve $(1-C_{w^e})\mu^e=F^e$.
Then $$\Norm2{\Sigma^e}{\mu^e}\le
c\Norm2\rb{(1-C_{w^e})^{-1}}\Norm2{\Sigma^e}{\hat F}.$$ But again
on $\rb$ we have $(1-\cht)(\mu^e\restriction\rb)=\hat F$ i.e.
$\hat\mu=\mu^e\restriction\rb$ solves $ (1-\cht)\hat\mu=\hat F$
and
$\Norm2{\rb}{\hat\mu}\le\Norm2{\Sigma^e}{\mu^e}\le\Norm2\rb{(1-C_{w^e})^{-1}}\Norm2{\Sigma^e}{\hat
F}$. Also it is clear that $\hat\mu$ is the unique solution of
$(1-C_{\hat w_\theta})\hat\mu=\hat F$.  For if
$(1-\cht)\tilde\mu=0$, $\tilde\mu\in L^2(\rb)$, then
$(1-C_{w^e})\tilde\mu^e=0$ where $\tilde\mu^e\in L^2(\Sigma^e)$ is
given by $\tilde\mu$ on $\rb$ and by $C(\tilde\mu(\hat
w_\theta^++\hat w_\theta^-))$ on $\Sigma^e\setminus\rb$. But then
as $(1-C_{w^e})^{-1}$ exists, $\tilde\mu^e=0$ and hence
$\tilde\mu=0$. It follows that we must prove (2.52) on $\Sigma^e$.

In view of the preceding remarks we must show that the RHP
$(\Sigma^e,v^e)$ localizes near $z_0$ and by the steepest descent
method introduced in [DZ1], this is done by taking into account
the signature table of $\re i\theta=\re
i(tz_0^2-t(z-z_0)^2)=t\im(z-z_0)^2$

\bigskip\bigskip
\font\thinlinefont=cmr5

\centerline{\beginpicture \setcoordinatesystem units
<.7000cm,.7000cm>
\setshadesymbol ({\thinlinefont .}) \setlinear
%
%
\linethickness= 0.500pt \setplotsymbol ({\thinlinefont .})
\putrule from  2.540 21.590 to  7.620 21.590
%
%
\linethickness= 0.500pt \setplotsymbol ({\thinlinefont .})
\putrule from  5.080 24.130 to  5.080 19.050
%
%
\put{$z_0$} [lB] at  5.239 21.114
%
%
\put{$\re i\theta<0$} [lB] at  2.651 22.860
%
%
\put{$\re i\theta>0$} [lB] at  5.7 22.860
%
%
\put{$\re i\theta>0$} [lB] at  2.651 20.161
%
%
\put{$\re i\theta<0$} [lB] at  5.7 20.161 \linethickness=0pt
\putrectangle corners at  2.515 24.155 and  7.645 19.025
\endpicture}

\centerline{Figure 2.53. Signature table for $\text{
 Re } i\theta$}
\bigskip
\noindent Guided by this table we must move the factors $e^{\pm
i\theta}$ into regions of the complex plane where they are
exponentially decreasing. It is here that IRHP2 plays a critical
role: rather than working with $(1-C_{w^e})^{-1}$ explicitly, we
use the analyticity inherent in the specification of IRHP2 to
perform the desired deformations of the problem. Again by
Proposition 2.14, in order to prove (2.52) it is enough to show
that the IRHP2 $$M_+^e=M_-^ev^e+F^e, \ \ M_{\pm}^e\in\dc2,$$ is
solvable for all $\hat F\in L^2$ and $$\Norm2{\Sigma^e}{M^e}\le
c\Norm2{\Sigma^e}{\hat F}\tag2.54$$ for all $t>0$ sufficiently
large and for all $x\in\rb$.

Define $\Phi_+,\Phi_-$ respectively on $\Sigma^e$ as follows:

$$\matrix\pmtwo10{[\bar
r]\delta^{-2}e^{-i\theta}}1,&\pmtwo1{[r]\delta^{2}e^{i\theta}}01,\
\ \ \ \ &\text{on} &z_0+\rb_+\\
\pmtwo1001,&\pmtwo10{[\bar r]\delta^{-2}e^{-i\theta}}1,\ \ \ \ \
&\text{on}
&z_0+e^{i\pi/4}\rb_+\\
\pmtwo1001,&\pmtwo1{-\left[\frac
r{1-|r|^2}\right]\delta^{2}e^{i\theta}}01,\ \ \ \ \ &\text{on}
&z_0+e^{-i\pi/4}\rb_-\\
\pmtwo1{-\left[\frac r{1-|r|^2}\right]\delta_+^{2}e^{i\theta}}01,
&\pmtwo10{-[\frac{\bar r}{1-|r|^2}]\delta_-^{-2}e^{-i\theta}}1, \
\ \ \ \ &\text{on}
&z_0+\rb_-\\
\pmtwo10{-[\frac{\bar
r}{1-|r|^2}]\delta^{-2}e^{-i\theta}}1,&\pmtwo1001,\ \ \ \ \
&\text{on}
&z_0+e^{i\pi/4}\rb_-\\
\pmtwo1{[r]\delta^{2}e^{i\theta}}01,&\pmtwo1001,\ \ \ \ \
&\text{on} &z_0+e^{-i\pi/4}\rb_+
\endmatrix$$and set $m_\pm^0=m_\pm^e\Phi_\pm$. Using the
properties of $r,\delta$ and the signature table for $\re
i\theta$, observe that $\Phi_\pm$ and $\Phi_\pm^{-1}$ are
uniformly bounded for all $x\in\rb,t\ge 0$. A direct calculation
shows that $m_\pm^0$ satisfy the relations $$m_+^0=m_-^0v^0+F^0$$
on $\Sigma^e$ where $F^0=F^e\Phi_+$ and
$$v^0=\cases\pmtwo1001,\ \ \ \ \ &\text{on }
\rb\\
\pmtwo10{-[\bar r]\delta^{-2}e^{-i\theta}}1,\ \ \ \ \ &\text{on }
z_0+e^{i\pi/4}\rb_+\\
\pmtwo1{\left[\frac r{1-|r|^2}\right]\delta^{2}e^{i\theta}}01,\ \
\ \ \ &\text{on }
z_0+e^{-i\pi/4}\rb_-\\
\pmtwo10{-[\frac{\bar r}{1-|r|^2}]\delta^{-2}e^{-i\theta}}1,\ \ \
\ \ &\text{on }
z_0+e^{i\pi/4}\rb_-\\
\pmtwo1{[r]\delta^{2}e^{i\theta}}01,\ \ \ \ \ &\text{on }
z_0+e^{-i\pi/4}\rb_+
\endcases\tag2.55$$
Moreover, $m_\pm^0\in\dc2$. To verify this, observe that
$\Phi_\pm$ are the boundary values of a bounded, analytic function
$\Phi$, say, in $\cb\setminus\Sigma^e$. For example, in the sector
${0<\arg z<\pi/4}$,  $\Phi(z)=\pmtwo10{[\bar
r]\delta^{-2}e^{-i\theta}}1$, etc. Writing
$m_\pm^e=C_{\Sigma^e}^\pm h$ for some $h\in L^2(\Sigma^e)$, we
choose $h^\ep\in L^2(\Sigma^e)$ such that $h^\ep\to h$ in
$L^2(\Sigma^e)$ as $\ep\downarrow0$, and set
$m^\ep(z)=(C_{\Sigma^e}h^\ep)(z),\  z\in\cb\setminus\Sigma^e$.
Applying Cauchy's theorem to $m^\ep(z)\Phi(z)$, we find
$m^\ep(z)\Phi(z)=(C_{\Sigma^e}(m_+^\ep\Phi_+-m_-^\ep\Phi_-))(z),\
z\in\cb\setminus\Sigma^e$. Evaluating this relation on $\Sigma^e$,
and then letting $\ep\downarrow0$, we obtain
$$m_\pm^0=m^\ep_\pm(z)\Phi_\pm(z)=(C^\pm_{\Sigma^e}(m_+^\ep\Phi_+-m_-^\ep\Phi_-))(z)
\in\dc2$$as desired. As $\Phi_\pm$ and $\Phi_\pm^{-1}$ are
uniformly bounded in $x$ and $t$, we see that (2.54) is equivalent
to showing that
the IRHP2
$$
M^0_+ = M^0_-v^0 + F^0,\quad M^0_\pm \in \dc2,
$$
is solvable for all $F^0\in L^2$ and
$$
\|M^0_\pm\|_{L^2(\Sigma^e)} \le c\|F^0\|_{\Sigma^e} \tag 2.56$$
for all $t>0$ sufficiently large and for all $x\in \rb$. The key
point to note is that the IRHP2 for $v^e$ has been deformed to a
problem with jump matrix $v^0$ where all the exponential factors
$e^{\pm i\theta}$ are now exponentially decreasing.

We need to prove some estimates on $\delta(z)$. Let
$\|r\|_{L^\infty(\rb)}\le\rho < 1$. From Proposition 2.38, we have
$\delta(z) = e^{(C_{\rb_-+z_0} \log 1-|r|^2)(z)}$, $z\in
\cb\setminus(\rb_-+z_0)$. Now where $\chi^0(s)$ denotes the
characteristics function of the interval $(z_0-1,z_0)$ we have
$z\in \cb\setminus(\rb_-+z_0)$
$$\align
(C_{\rb_-+z_0} \log(1-|r|^2))(z) &= \int^{z_0}_{-\infty} \{\log(1-|r(z)|^2) - \log(1-|r(z_0)|^2) \chi^0(s)(s-z_0+1)\} \frac{ds}{2\pi i(s-z)}\\
&\quad + \log(1-|r(z_0)|^2) \int^{z_0}_{z_0-1} \frac{s-z_0+1}{s-z} \frac{ds}{2\pi i}\\
&= \beta(z,z_0) + i\nu(z_0)(1-(z-z_0 + 1) \log(z-z_0+1)\\
&\quad + (z-z_0) \log(z-z_0) + \log(z-z_0))
\endalign$$
where $\nu(z_0) = -\frac1{2\pi} \log(1-|r(z_0)|^2)$ as in (1.2),
and
$$
\beta(z,z_0) = \int^{z_0}_{-\infty} \{\log(1-|r(z)|^2) -
\log(1-|r(z_0)|^2) \chi^0(s) (s-z_0+1)\} \frac{ds}{2\pi i(s-z)}.
$$
Clearly $\beta_\pm(\cdot,z_0)$ lies in the Sobolev space $H^1(\rb)
= H^{1,0}$, and a direct calculation shows that
$$
\|\beta_\pm\|_{H^1(\rb)} \le \frac{c\|r\|_{H^{1,0}}}{1-\rho}.
$$
Focusing on the ray $L_{-\pi/4} = z_0 + e^{-i\pi/4} \rb_+ =\{z =
z_0 + ue^{-i\pi/4}, u\ge 0\}$, we find similarly

$$
\|\beta\|_{H^1(L_{-\pi/4})} \le \frac{c\|r\|_{H^{1,0}}}{1-\rho}.
$$
It follows then by standard Sobolev estimates that on $L_{-\pi/4}$
$$
\left\{\matrix
\text{(i)}&\beta(z,z_0) \text{ is continuous up to } z=z_0\hfill\\
\text{(ii)}&\|\beta(\lozenge,z_0)\|_{L^\infty(L_{-\pi/4})} \le c\|r\|_{H^{1,0}}/1-\rho\hfill\\
\text{(iii)}&|\beta(z,z_0) - \beta(z_0,z_0)| \le {\dsize\frac{c\|r\|_{H^{1,0}}}{1-\rho}} |z-z_0|^{1/2}\hfill
\endmatrix\right.\tag 2.57
$$
where the constant $c$ is independent of $x\in\rb$ and $t>0$.
Write
$$
\delta = \delta_0\delta_1\tag 2.58$$
where
$$
\delta_0 = e^{\beta(z_0,z_0) + i\nu(z_0)} (z-z_0)^{i\nu(z_0)},\ \
\ \delta_1 = e^{\zeta(z,z_0)}$${and} $$\zeta(z,z_0) =
(\beta(z,z_0) - \beta(z_0,z_0)) + i\nu(z_0) ((z-z_0) \log(z-z_0) -
(z-z_0+1) \log(z-z_0+1)). $$
 We want to estimate
$$\align
D(z) &\equiv [r]\delta^2e^{i\theta} - r(z_0) \delta^2_0e^{i\theta}\tag 2.59\\
&= \delta^2_0 r(z_0) e^{i\theta} \left(\frac{\delta^2_1}{(1+i(z-z_0))^2} - 1\right)
\endalign$$
on $L_{-\pi/4}$ (cf.\ (2.55)). Write $D(z) = D_1(z) + D_2(z)$ where
$$\align
D_1(z) &= \delta^2 r(z_0) e^{itz^2_0} e^{-it(z-z_0)^2} \frac{(-2i(z-z_0) + (z-z_0)^2)}{(1+i(z-z_0))^2}\\
\intertext{and}
D_2(z) &= \delta^2_0 r(z_0) e^{itz^2_0} e^{-it(z-z_0)^2} (\delta^2_1(z)-1).
\endalign$$
Now as $\beta(z_0,z_0) \in i\rb$, $|\delta_0| = |e^{i\nu
\log(z-z_0)}| = e^{\nu(z_0)\pi/4}$ and so as $|\delta(z)| \le
(1-\rho)^{-1/2}$, we learn that $|\delta_1(z)| \le
\frac{e^{-\nu(z_0)\pi/4}}{(1-\rho)^{1/2}} \le
\frac1{(1-\rho)^{1/2}}$. On $L_{-\pi/4} = \{z = z_0 +
ue^{-i\pi/4}$, $u\ge 0\}$, $|D_1(z)| \le \frac\rho{1-\rho}
e^{-tu^2}(2u+u^2) \le \frac{c}{1-\rho} \frac{e^{-t/2
u^2}}{t^{1/2}}$ for $t\ge 1$, where $c = \sup\limits_{w>0}
e^{-w^2/2}(2w+w^2)$. Also
$\left\|\frac{d\xi}{dz}\right\|_{L^2(L_{\pi/4})} \le c$ where $c$
is independent of $x$ and $t\ge 0$, by the $H^1$ property of
$\beta(z,z_0)$ and the properties of the logarithm:\ hence
$|\zeta(z,z_0)| \le c|z-z_0|^{1/2}$, $z\in L_{-\pi/4}$. It follows
that for $z\in L_{-\pi/4}$,
$$\align
|\delta^2_1(z) - 1| &= \left|\int^1_0 \frac{d}{ds} e^{2s\zeta(z,z_0)}ds\right| \le c|z-z_0|^{\frac12} \max_{0\le s\le 1} |e^{2\zeta(z,z_0)}|^s\\
&\le c|z-z_0|^{1/2} \max_{0\le s\le 1}|\delta_1(z)|^{2s} \le \frac{c}{1-\rho} |z-z_0|^{1/2},
\endalign$$
and we conclude that on $L_{-\pi/4}$, $|D_2(z)| \le ce^{\nu\pi/2} \frac\rho{1-\rho} e^{-tu^2} u^{1/2}  \le \frac{ce^{-t/2 u^2}}{t^{1/4}}$. Assembling the above estimates we have for $t\ge 1$ and $z\in L_{-\pi/4}$,
$$
|D(z)| \le \frac{ce^{(-t/2)|z-z_0|^2}}{t^{1/4}}\tag 2.60$$ The
situation on the other rays in $\Sigma^e\backslash \rb$ is
similar. Set
$$
\left\{\matrix
v^{\#} = I\quad \text{on}\quad \rb\hfill\\
v^{\#} = \left(\matrix 1&0\\ -\bar r(z_0)\delta^{-2}_0 e^{-i\theta}&1\endmatrix\right) \quad \text{on}\quad z_0 + e^{i\pi/4} \rb_+\hfill\\
v^{\#} = \left(\matrix 1&(r(z_0)/(1-|r(z_0)|^2))\delta^2_0e^{i\theta}\\ 0&1\endmatrix\right) \quad \text{on}\quad z_0 + e^{-i\pi/4} \rb_-\hfill\\
v^{\#} = \left(\matrix 1&0\\ -(\ovl{r(z_0)}/(1- |r(z_0)|^2))
\delta^{-2}_0 e^{-i\theta}&1\endmatrix \right) \quad
\text{on}\quad z_0 + e^{i\pi/4} \rb_-
\hfill\\
v^{\#} = \left(\matrix 1&r(z_0) \delta^2_0 e^{i\theta}\\
0&1\endmatrix\right)\quad \text{on}\quad z_0 + e^{-i\pi/4}
\rb_+.\hfill
\endmatrix\right.\tag 2.61
$$
We have proved the following.

\proclaim{Lemma 2.62}
For $1\le p \le \infty$,
$$
\|v^{\#} - v^0\|_{L^p(\Sigma^e)} \le
\frac{c}{t^{\frac14+\frac1{2p}}}\tag 2.63$$ uniformly for $t\ge 1$
and all $x\in\rb$.
\endproclaim
\n {\bf Remarks}
\roster
\item If we knew that $r$ had more smoothness, say $r''\in L^2$, then we could replace estimate $|\beta(z,z_0) - \beta(z_0,z_0)| \le c|z-z_0|^{1/2}$ by $|\beta(z,z_0) - \beta(z_0,z_0)| \le c|z-z_0|$, which would lead eventually to the bound $\frac{\log t}{t^{1/2+1/2p}}$ in (2.63).
\item If $r$ had more decay, say $r\in H^{1,1}_1$, then one can easily show that $\beta(z_0,z_0) + i\nu(z_0) = -\frac1{2\pi i} \int^{z_0}_{-\infty} \log(z_0-s) d \log 1-|r(s)|^2$. We will need this fact later.
\endroster

Set $w^{\#} = (0, w^{\#+} = v^{\#} -I), w^0 = (0,w^{0+} = v^0-I)$ on $\Sigma^e$. We will show shortly that $(1-C_{w^{\#}})^{-1}$ exists on $L^2(\Sigma^e)$ for all $x,t$ and
$$
\|(1-C_{w^{\#}})^{-1}\|_{L^2(\Sigma^e)} \le c\tag 2.64$$
where $c$ is independent of space and time. By (2.63), for $t\ge 1$,
$$
\|C_{w^0} - C_{w^{\#}}\|_{L^2}  = \|C^-\loz (v^0-v^{\#})\|_{L^2}
\le c\|v^0-v^{\#}\|_{L^\infty} \le \frac{c}{t^{1/4}},
$$
and we conclude by the resolvent identity that $(1-C_{w^0})^{-1}$ exists in $L^2(\Sigma^e)$ for $t$ sufficiently large and that
$$
\|(1-C_{w^0})^{-1}\|_{L^2(\Sigma^e)} \le c,\quad
\|(1-C_{w^0})^{-1} - (1-C_{w^{\#}})^{-1}\|_{L^2} \le
\frac{c}{t^{1/4}}\tag 2.65$$ uniformly for all $x\in \rb$ and
$t\ge t_0$, say. But then again by Proposition (2.14), and the
Corollary to Proposition (2.14), this proves (2.56) and hence,
retracing the argument, we obtain the desired bound in (2.50).

Now in order to prove (2.64), observe that if we retrace the steps
for $v^\#$ rather than $v^0$, then instead of obtaining $\hat
v_\theta$ on $\rb$ we would obtain $\hat v{}^\#_0$, where
$$\align
\hat v{}^\#_\theta &= \left(\matrix 1&r(z_0) \delta^2_0 e^{i\theta}\\ 0&1\endmatrix\right) \left(\matrix 1&0\\ -\ovl{r(z_0)} \delta^{-2}_0 e^{-i\theta}&1\endmatrix\right),\qquad z>z_0,\\
&= \left(\matrix 1&0\\ -\left(\frac{\bar r(z_0)}{1-|r(z_0)|^2}\right) \delta^{-2}_{0-} e^{-i\theta}&1\endmatrix\right) \left(\matrix 1&\frac{r(z_0)}{1-|r(z_0)|^2} \delta^2_{0+} e^{i\theta}\\ 0&1\endmatrix\right),\qquad z<z_0,
\endalign$$
(cf.\ (2.48)). Multiplying out the matrices, we find
$$
\hat v{}^\#_\theta(z) = \delta^{\sigma_z}_{0-} \left(\matrix
1-|r(z_0)|^2&r(z_0)e^{i\theta(z)}\\ -\ovl{r(z_0)}
e^{-i\theta(z)}&1\endmatrix\right) \delta^{-\sigma_3}_{0+},\qquad
z\in\rb,\tag 2.66$$ where again $\sigma_3 = \left(\matrix 1&0\\
0&-1\endmatrix\right)$. Set $\widehat w{}^\#_0 = (0, \widehat
w{}^{\#+}_\theta = \hat v{}^\#_\theta - I)$. By Proposition
(2.14), and its Corollary, in order to complete the proof of
(2.64), and hence the proof of Proposition (2.49), it is enough to
show that $(1-C_{\widehat w{}^\#_\theta})^{-1}$ exists in
$L^2(\rb)$ for all $x,t$ and
$$
\|(1-C_{\widehat w{}^\#_\theta})^{-1}\|_{L^2} \le c \tag2.67$$
uniformly for all $x,t$. But by the analyticity and boundedness
properties of $\delta_0 = \delta_0(z)$, it follows that the IRHP2
for $\hat v{}^\#_\theta$ is equivalent to the IRHP2 for $v_1 =
\left(\matrix 1-|r_1|^2&r_1\\ -\bar r_1&1\endmatrix\right)$ where
$r_1(z) = r(z_0) e^{i\theta(z)}$ (cf.\ with the situation above
for $v^e$ and $v^0$). But as $\|r_1\|_{L^\infty} = |r(z_0)| \le
\|r\|_{L^\infty} < 1$, it follows from (2.30) that for $w_1 =
(w^-_1,w^+_1) = \left(\left(\matrix 0&r_1\\ 0&0\endmatrix\right),
\left(\matrix 0&0\\ -\bar r_1&0\endmatrix\right)\right),
(1-C_{w_1})^{-1}$ exists in $L^2$ and
$$
\|(1-C_{w_1})^{-1}\|_{L^2} \le \frac1{1-\|r\|_{L^\infty}}
\tag2.68$$
for all $x,t$. The bound (2.67) then follows by two applications of Proposition 2.14. This completes the proof of Proposition (2.49).

\remark{Remark \rom{2.69}}
The above proof of (2.50) should be contrasted with the proof of the same fact in \cite{DZ2}, which relies in turn on various intricate formulae taken from \cite{DZ1}.
\endremark

\head
3. Scattering and Inverse Scattering for $\bold{q\in H^{1,1}}$
\endhead

The scattering/inverse scattering theory for the ZS--AKNS operator
$T  = \partial_x - \left(iz\sigma + \left(\matrix 0&q(x)\\
\ovl{q(x)}&0\endmatrix\right)\right)$ is well known. An account of
the classical theory is contained in \cite{ZS}, \cite{AKNS}, and
the Riemann--Hilbert point of view can be inferred, as a special
case, from the general theory in \cite{BC}.

In this section, we present a somewhat complete, mostly
self-contained, sketch of the general theory for the operator $T$.
Our main focus is on a detailed analysis (following \cite{Z1}) of
the mapping properties of the scattering map $\cl R$, showing, in
particular, that $\cl R$ is a bijection from $H^{1,1}$ onto
$H^{1,1}_1$. Enough of the general theory is described, however,
so that the interested reader should be able to fill in the
missing details without too much effort.

The principal objects of study in the theory are eigensolutions $\psi = \psi(x,z)$ of the ZS--AKNS operator
$$
\left(\partial_x - \left(iz\sigma + \left(\matrix  0&q(x)\\ \bar q(x)&0\endmatrix \right)\right)\right) \psi = 0.\tag3.1$$
Setting
$$
m = \psi e^{-ixz\sigma},\tag3.2$$
equation (3.1) takes the form
$$
\partial_xm = iz \text{ ad } \sigma(m) + Qm,\qquad Q = \left(\matrix  0&q\\ \bar q&0\endmatrix \right),
\tag3.3$$
where ad $A(B) = [A,B] = AB-BA$. Under exponentiation we have $e^{\text{ad } A} B = \sum^\infty_{n=0} \frac{(\text{ad } A)^n(B)}{n!} = e^ABe^{-A}$ (cf.\ $A_\theta$ following (2.35)). The theory of ZS--AKNS \cite{ZS}, \cite{AKNS} is based on the following two Volterra integral equations for real $z$,
$$
m^{(\pm)} (x,z) = I + \int^x_{\pm\infty} e^{i(x-y)z\text{ ad }
\sigma} Q(y) m^{(\pm)} (y,z)dy \equiv I +
K_{q,z,\pm}m^{(\pm)}.\tag3.4$$ By iteration, one sees that these
equations have bounded solutions continuous for both $x$ and real
$z$ when $q\in L^1(\rb)$. The matrices $m^{(\pm)}(x,z)$ are the
unique solutions of (3.4) normalized to the identity as
$x\to\pm\infty$. The following are some relevant results of
ZS--AKNS theory:
\medskip

\n (3.5a)\quad
There is a continuous matrix function $A(z)$ for real $z$, $\det A(z) = 1$, defined by $\psi^{(+)} = \psi^{(-)}A(z)$, where $\psi^{(\pm)} = m^{(\pm)} e^{ixz\sigma}$ and $A$ has the form $A(z) = \left(\matrix  a&\bar b\\ b&\bar a\endmatrix \right)$.

\n (3.5b)\quad $a$ is the boundary value of an analytic function,
also denoted by $a$, in the upper half-plane $\cb_+$:\ $a$ is
continuous and non-vanishing in $\ovl{\cb_+}$, and
$\lim\limits_{z\to \infty} a(z)=1$.
$$
\left\{\matrix
a(z) = \det(m^{(+)}_1, m^{(-)}_2) = 1 - {\ds \int_{\rb}} q(y) m^{(+)}_{21}(y,z)dy = 1 + {\ds\int_{\rb}} \ovl{q(y)} m^{(-)}_{12}(y,z)dy,\hfill\\
b(z) = e^{ixz} \det(m^{(-)}_1,m^{(+)}_1) = -{\ds\int_{\rb}}
\ovl{q(y)} e^{iyz} m^{(+)}_{11} (y,z) = -{\ds\int_{\rb}}
\ovl{q(y)} e^{iyz} m^{(-)}_{11}(y,z)dy,\hfill
\endmatrix\right.
\tag3.5c$$
where $m^{(\pm)} = (m^{(\pm)}_1, m^{(\pm)}_2) = \left(\matrix  m^{(\pm)}_{11}&m^{(\pm)}_{12}\\ m^{(\pm)}_{21}&m^{(\pm)}_{22}\endmatrix \right)$.

\n (3.5d)\quad The {\it reflection coefficient\/} $r$ is defined
by $-\bar b/\bar a$. As $\det A=1$, $|a|^2 - |b|^2 = 1$ so that
$|a|\ge 1$ and $|r|^2=1-|a|^{-2}<1$. Together with (3.5b), this
implies $\|r\|_{L^\infty(\rb)} <1$. The {\it transmission
coefficient\/} $t(z)$ is defined by $1/a(z)$. Thus
$\|t\|_{L^\infty(\rb)}\le 1$ and $|r(z)|^2 + |t(z)|^2 =1$.

\remark{Remark 3.6}
It is easy to verify that if we set $z_0 = +\infty$ in (2.39), then $\delta(z) = t(z)$.
\endremark

Now $m^{(+)}_1(x,z)$ and $m^{(-)}_2(x,s)$, the first and second
columns of $m^{(+)}$ and $m^{(-)}$ respectively, have analytic
continuations to $\cb_+$. Moreover if we set
$$
m = m(x,z) = \left(\frac{m^{(+)}_1(x,z)}{a(z)},
m^{(-)}_2(x,z)\right),\qquad z\in \cb_+,
$$
then $m(x,z)$ is analytic in $\cb_+$, continuous in $\ovl{\cb_+}$
and
$$
\left\{\matrix
m(x,z) \to I\quad \text{as}\quad x\to-\infty\hfill\\
m(x,z) \quad \text{is bounded as}\quad x\to +\infty.\hfill
\endmatrix\right.\tag3.7$$
Using the fact that the substitution $m(x,z) \to \left(\matrix  0&1\\ 1&0\endmatrix \right) \ovl{m(x,\bar z)} \left(\matrix  0&1\\ 1&0\endmatrix \right)$ preserves solutions of (3.3), we define
$$
m(x,z) \equiv \left(\matrix  0&1\\ 1&0\endmatrix \right) \ovl{m(x,\bar z)} \left(\matrix  0&1\\ 1&0\endmatrix \right)
$$
for $z\in \cb_-$, and by continuity in $\ovl{\cb_-}$. With these
definitions $\psi(x,z) = m(x,z)e^{ixz\sigma}$, $z\in
\cb\backslash\rb$, is the (unique) solution of (3.1) of
\cite{BC}--type mentioned in the Introduction. Orienting $\rb$
from $-\infty$ to $+\infty$, we let $\psi_\pm(x,z)$, $m_\pm(x,z)$
denote the boundary values of $\psi(x,z) = m(x,z)e^{ixz\sigma}$,
$m(x,z)$ from $\cb_\pm$ respectively. By uniqueness of solutions
of (3.1), we must have $\psi_+(x,z) = \psi_-(x,z)v$ for some
$v=v(z)$ independent of $x$. Direct calculation shows that $v(z) =
\left(\matrix  1-|r(z)|^2&r(z)\\ -\ovl{r(z)}&1\endmatrix \right)$
where $r=-\bar b/\bar a$ is the reflection coefficient in (3.5d).
For fixed $x,m(x,z)$ solves the normalized RHP (1.4) in the
following precise sense
$$
\left\{\matrix \bullet\hfill&m(x,z)  \text{ is analytic in
$\cb\backslash \rb$ and continuous up to the boundary}
\hfill\\
\bullet\hfill&m_+(x,z) = m_-(x,z)v_x(z) \text{ for $z\in \rb$, where } v_x = e^{ix \text{ ad } \sigma} v\hfill\\
\bullet\hfill&m(x,z)\to I \text{ uniformly as $z\to \infty$ in }
\ovl{\cb_+} \cup \ovl{\cb_-}.\hfill
\endmatrix\right.\tag3.8$$
The above considerations require only that $q\in L^1(\rb)$. We are
interested in the case that $q\in H^{0,1}\subset L^1$, and for
such $q$ we will see $m_\pm$ also solves the normalized RHP in the
sense (2.9),
$$
m_+ = m_-v_x\quad \text{on}\quad \rb,\qquad m_\pm -I\in\dc2.
$$
In fact, even more is true, viz.\ $m_\pm -I\in \partial C(H^1)$. In other words $m_\pm = C^\pm h$, where $h$ lies in the Sobolev space $H^1 = H^{1,0}$.

As noted in the Introduction, in \cite{Z1}, the author proved that $q\mapsto r=\cl R(q)$ is a bi-Lipschitz map from $H^{k,j}$ to $H^{j,k}_1$ for $k\ge 0$, $j\ge 1$. We now prove this fact in the special case $k=j=1$, which is, of course, the main case of interest in this paper. For simplicity, we will only prove that the map is bijective. Our method is based on the analysis of certain rational forms of some linear operators which depend linearly on $q,\bar q$ or $r,\bar r$ and it will be clear to the reader that the maps $\cl R$ and $\cl R^{-1}$ are indeed Lipschitz.

The bijectivity of $\cl R$ is covered by the following four
theorems. In the following if $M$ is a measure space and $B$ is a
Banach space, then $B\otimes L^p(M) \equiv L^p(M\to B)$ denotes
the space of $B$-valued $L^p$ functions with norm $\|f\|_{B\times
L^p(M)} = \left\|\,\|f\|_B\,\right\|_{L^p(M)}$.

\proclaim{Theorem 3.9}
$\cl R$ maps $H^{0,1}$ into $H^{1,0}$.
\endproclaim

\demo{Proof}
By the Wronskian formula (3.5c) and the fact that $|a(z)|\ge 1$, we only need to prove that $m^{(\mp)} (x,\cdot) - I\in H^{1,0}$ for any fixed $x$, say $x=0$. We will provide the proof only for $m^{(+)}$. Clearly, $\|K_{q,z,+}\|_{L^\infty(dx)\to L^\infty(dx)}$ $\le \|Q\|_{L^1}$. The operator $K_{q,z,+}$ can also be lifted to a bounded operator $K_{q,+}$ from $L^2(dz) \otimes L^\infty(dx)$ to $L^2(dz)\otimes L^\infty(dx)$:
$$
\|K_{q,+}f\|_{L^2(dz)\otimes L^\infty(dx)}\le \|Q\|_{L^1}
\|f\|_{L^2(dz)\otimes L^\infty(dx)}.
$$
The standard iteration method for Volterra integral equation gives the estimate
$$
\|(1-K_{q,+})^{-1}\|_{L^2(dz)\otimes L^\infty(dx)\to L^2(dz) \otimes L^\infty(dx)}\le e^{\|Q\|_{L^1}}.
$$

By Fourier theory and Hardy's inequality,
$$\align
\|K_{q,+}I\|_{L^2(dz)\otimes L^\infty(dx)} &= \sqrt{2\pi} \left\|\left(\int^{\langle x\rangle}_{+\infty} |Q|^2dy\right)^{\frac12} \right\|_{L^\infty(dx)} \le \sqrt{2\pi} \|Q\|_{L^2}\\
\|K_{q,+}I\|_{L^2(dz) \otimes L^2_{\rb_+}(dx)} &= \sqrt{2\pi}
\left\|\left(\int^{\langle x\rangle}_{+\infty}
|Q|^2dy\right)^{\frac12} \right\|_{L^2_{\rb_+}(dx)} \le
\sqrt{2\pi} \|Q\|_{H^{0,\frac12}} \le \sqrt{2\pi}\|Q\|_{H^{0,1}}.
\endalign$$

Thus
$$
\|m^{(+)} - I\|_{L^2(dz)\otimes L^\infty(dx)} =
\|(1-K_{q,+})^{-1}K_{q,+}I\|_{L^2(dz)\otimes L^\infty(dx)} \le
ce^{\|Q\|_{L^1}} \|Q\|_{L^2}.\tag3.10$$ Estimating the RHS of
(3.4), and using Hardy's inequality \cite{HLP}, we see that in
fact
$$
m^{(+)}  -I\in L^2(dz) \otimes (L^\infty\cap L^2_{\rb_+}(dx)).
\tag3.11$$

Since $\partial_zm^{(+)}$ satisfies the equation
$$\align
\partial _z m^{(+)}(x,z) &= \int^x_{+\infty} i(x-y) \text{ ad } \sigma e^{i(x-y)z \text{ ad } \sigma} Q(y) m^{(+)}(y,z)dy\\
&\quad + \int^x_{+\infty} e^{i(x-y)z \text{ ad } \sigma} Q(y) \partial _z m^{(+)} (y,z)dy,
\endalign$$
the function $\widehat m {\buildrel \text{def}\over =} (\partial_z-ix \text{ ad } \sigma)m^{(+)}$ satisfies the equation
$$\align
\widehat m(x,z) &= -i \int^x_{+\infty} ye^{i(x-y)z \text{ ad } \sigma} (\text{ad } \sigma Q(y))m^{(+)} (y,z)dy\\
&\quad + \int^x_{+\infty} e^{i(x-y)z \text{ ad } \sigma} Q(y) \widehat m(y,z)dy\\
&\equiv h_1 + h_2 + K_{q,+}\widehat m,
\endalign$$
where
$$
h_1 = -i \int^{\langle x\rangle}_{+\infty} ye^{i(\langle x\rangle
-y)\langle z\rangle \text{ ad } \sigma} (\text{ad } \sigma Q(y))dy
\in L^2(dz) \otimes L^\infty(dx),
$$
and
$$
h_2 = -i \int^{\langle x\rangle}_{+\infty} ye^{i(\langle x\rangle -y)\langle z\rangle \text{ ad } \sigma} (\text{ad } \sigma Q(y)) (m^{(+)}(y, \langle z\rangle) -I)dy \in L^2(dz) \otimes L^\infty(dx)
$$
by (3.11) and the fact that $Q\in H^{0,1}$. Hence $h_1+h_2\in
L^2(dz)\otimes L^\infty(dx)$ and
$$
\widehat m = (1-K_{q,+})^{-1} (h_1+h_2) \in L^2(dz) \otimes L^\infty(dx).
$$
In particular $\partial_zm^{(+)}(x=0,z) = \partial_x\widehat m(0,z)\in L^2(dz)$ and by (3.10), $m^{(+)}(0,z)-I\in L^2(dz)$. Analogous estimates for $m^{(-)}$ show that $\partial_zm^{(-)}(0,z)$ and $m^{(-1)}(0,z)-I$ lie in $L^2(dz)$, and hence $r\in H^{1,0}$ by (3.5bc).
\enddemo

Note that the above calculations show that $m^{(+)}(x,z)-I$ and
$\partial_z m^{(+)}(x,z)\in L^2(dz) \otimes L^\infty(dx)$. There
is a similar calculation for $m^{(-)}(x,z)$. It is then easy to
show that for each $x\in \rb$, $m_\pm(x,\cdot)-I\in \partial
C(H^1)$, as advertised above.

We are now ready to consider the regularity of $q$.

\proclaim{Theorem 3.12}
If $q\in H^{1,1}$, then $r\in H^{1,1}_1$.
\endproclaim

\demo{Proof}
Since $r\in H^{1,0}$, we only need to study its decay behavior as $|z|\to \infty$. As $|a|\ge 1$, we only need to study the decay behavior of $b$ at $\infty$. We use the formulae (3.5c),
$$
b = -\int \bar qe^{iyz} m^{(-)}_{11} = -\int \bar q e^{iyz} (m^{(-)}_{11}-1) - \int \bar q e^{iyz}.
$$
Since the second term on the RHS is clearly in $H^{1,1}$, we only need to show that $\int \bar q e^{iyz}(m^{(-)}-I)$ is in $H^{0,1}$.

Denote $K_{q,z,-}$ by $K$. Noting that $\text{ad } \sigma$ is invertible on off-diagonal matrices we integrate by parts to obtain
$$
(KI)(x) = -(iz \text{ ad } \sigma)^{-1} Q+(iz \text{ ad } \sigma)^{-1} \int^x_{-\infty} e^{i(x-y)z \text{ ad } \sigma} Q'(y)dy \equiv h_1 + h_2.
$$
Note that we are only interested in the decay at $z=\infty$; the singularity at $z=0$ is irrelevant.

We write
$$
m^{(-)} - I = (1-K)^{-1} KI = h_1 + Kh_1 + (1-K)^{-1} K^2h_1 + (1-K)^{-1} h_2\equiv h_1 + g_1 + g_2 +g_3.
$$
Thus we need to show that
$$
\int \bar q e^{iyz} (h_1+g_1+g_2+g_3)
$$
is in $H^{0,1}$. It is obvious that $\int \bar qe^{iyz} h_1$ is in $H^{0,2}$.

Since $Q \text{ ad } \sigma Q$ is diagonal,
$$
(K(z^{-1} \text{ ad } \sigma Q)(x) = z^{-1} \int^x_{-\infty} e^{i(x-y)z \text{ ad } \sigma} Q \text{ ad } \sigma Q(y)dy = z^{-1} \int^x_{-\infty} Q \text{ ad } \sigma Q\ dy.
$$
Now $\int \bar qe^{iyz}g_1$ is in $H^{0,2}$ because $\bar q \int^{\langle x\rangle}_{-\infty} Q \text{ ad } \sigma Q \in H^{1,0}(dx)$.

Since $Q \text{ ad } Q$ is diagonal, and $Q \int^x_{-\infty} Q \text{ ad } \sigma Q$ is off-diagonal,
$$
(K^2(z^{-1} \text{ ad } \sigma Q)(x) = z^{-1} \int^x_{-\infty} e^{i(x-y)z \text{ ad } \sigma} Q(y) \int^y_{-\infty} Q \text{ ad } \sigma Q(u)du
$$
is in $H^{0,1} \otimes L^\infty(dx)$ and so is $g_2$. We see that $\int \bar qe^{iyz} g_2$ is also in $H^{0,1}$.

But, $h_2$ is in $H^{0,1} \otimes L^\infty(dx)$ and so is $g_3$. Thus $\int \bar q e^{iyz} g_3$ is also in $H^{0,1}$.

Finally we indeed have $r\in H^{1,1}_1$. As $b,\partial_z b \in L^2(dz)$, $b(z) \to 0$ as $z\to \infty$ and hence $|a(z)|^2 = 1 + |b(z)|^2\to 1$. Clearly $|a(z)|\ge 1$ and hence the bound $\|r\|_{L^\infty(dz)} < 1$ follows from the formula $|r(z)|^2 = 1 - \frac1{|a(z)|^2}$.
\enddemo

We now construct the inverse scattering transform from $H^{1,1}_1$ to $H^{1,1}$. Let $r\in H^{1,1}_1$ be given.

The inverse problem is formulated as a normalized RH problem
$(\rb,v_x) = (\rb,v_x)_{L^2}$ where $v = v^{-1}_-v_+ =
\left(\matrix  1-|r|^2&r\\ -\bar r&1\endmatrix \right)$. This
problem is adapted for studying the decay behavior of $q$ as $x\to
-\infty$. To obtain an equivalent RH problem adapted for studying
the decay behavior of $q$ as $x\to\infty$, we set $\check{m} =
m\delta^{-\sigma_3}$, $\delta = \delta_{z_0=+\infty}= \exp C_{\rb}
\log(1-|r|^2)$, and $\check{v} = \delta^{\sigma_3}_-
v\delta^{-\sigma_3}_+$ \cite{BC}. In the following Lemma and
Theorem, we prove the decay behavior of $q$ as $x\to -\infty$. The
decay behavior of $q$ as $x\to\infty$ can be obtained in a similar
manner using the RH problem for $\tilde v$.

The normalized RHP $(\rb,v_x)$,
$$
m_+ = m_-v_x \quad \text{on}\quad \rb,\qquad m_\pm \in I  +
\partial C(L^2)
$$
is solvable for each $x$ by Proposition 2.31, and moreover, by
(2.32)--(2.33), $m_\pm = I+C^\pm(\mu(w^+_x + w^-_x))$ where $\mu
\in I + L^2$ is the unique solution of $(1-C_{w_x})\mu = I$.

But more is true:\ if $r\in H^{1,0}$, then $\mu-I\in H^1(dz)$.
Consider first the case where $r\in \cl C^\infty_0(\rb)$ and for
any $\vp\in \rb$ write $\mu^\vp(z) = \mu(z+\vp)$. Then the
equation for $\mu$ implies $(\mu^\vp-\mu)/\vp =
(1-C_{w^\vp_x})^{-1} C_{((w^\vp_x-w_x)/{\vp})}\mu$, and using the
fact that $\mu\in I+L^2$ and $\frac{w^\vp_x-w_x}\vp \to
\frac{d}{dz} w^\vp_x$ in $L^2\cap L^\infty(dz)$ as $\vp \downarrow
0$, we infer that $\mu'$ exists in $L^2(dz)$ and $\mu' =
(1-C_{w_x})^{-1} C_{w'_x}\mu$. But then as $\mu-I\in L^2$,
$\mu'\in L^2$ we have by Sobolev, $\|\mu-I\|_{L^\infty} \le
\gamma^{-1}\|\mu-I\|_{L^2} + \gamma\|\mu'\|_{L^2}$ for any
$\gamma>0$. Inserting this bound in $\|\mu'\|_{L^2} \le
\frac{c}{1-\|r\|_{L^\infty}} (I+
\|\mu-I\|_{L^\infty})\|r\|_{H^{1,0}}$, we obtain the {\it a
priori\/} bound, $\|\mu'\|_{L^2} \le
\frac{c\|r\|_{H^{1,0}}}{\left(1-\left(\frac{c\gamma
\|r\|_{H^{1,0}}}{1-\|r\|_{L^\infty}}\right)\right)}
\left(1+\frac1\gamma \|r\|_{H'}\right)$ for $\gamma$ sufficiently
small.  Note that the constant $c=c(x)$ grows at most linearly
with $x$ and therefore, for given $0<\rho<1$, $\lambda>0$,
$\gamma$ may be chosen uniformly for $\|r\|_{H^{1,0}}\le \lambda$,
$\|r\|_{L^\infty}\le\rho$ and all $x$ in compacta. Now given $r\in
H^{1,0}_1$, $\|r\|_{H^{1,0}} < \lambda$, $\|r\|_{L^\infty} < \rho
< 1$, we may choose $r_n\in \cl C^\infty_0(\rb)$,
$\|r_n\|_{H^{1,0}} \le \lambda$, $\|r_n\|_{L^\infty} \le \rho$
such that $r_n\to r$ in $H^{1,0}$. By the above calculations and
remarks it follows that $\mu-I\in H^1(dz)$ as desired. Moreover it
is clear that $\|\mu(x,\lozenge)-I\|_{H^1(dz)}$ is bounded for $x$
in compact sets.

Knowing that $\mu\in I+H^1(dz)\subset I+L^\infty(dz)$, we may now
differentiate $(1-C_{w_x})\mu=I$ with respect to $x$ to conclude
that $\frac{d}{dx}\mu = (1-C_{w_x})^{-1}(C_{(dw_x/dx)}\mu) \in
L^2(dz)$. But then
$$\frac{d}{dx}(M_\pm -I) = C^\pm \left[\left( \frac{d\mu}{dx}\right) (w^+_x+w^-_x) + \mu \left(\frac{d}{dx} (w^+_x+w^-_x)\right)\right] \in \partial C(L^2),$$
 and differentiating the jump relation $m_+  = m_-v_x$, we obtain $\frac{d}{dx} m_+ + iz[m_+,\sigma] = \left(\frac{d}{dx} m_- + iz[m_-,\sigma]\right)v_x$. A simple calculation shows that $iz[m_\pm,\sigma] = Q+C^\pm(i[\mu \lozenge(w^+_x+w^-_x),\sigma])$, where
$$
Q = Q(x) = \frac{\text{ad } \sigma}{2\pi} \int_{\rb} \mu(w^+_x +
w^-_x)dz.\tag3.13$$ It follows that $M_\pm \equiv \frac{d}{dx}
m_\pm -iz[\sigma, m_\pm] - Qm_\pm \in \partial C(L^2)$ and $M_+ =
M_-v_x$. But then $M_\pm=0$ by Proposition 2.14. Thus $m_\pm$
solve the differential equation $\frac{d}{dx} m_\pm = iz[\sigma,
m_\pm] + Qm_\pm$.

By the symmetry property of $v_x$, $v_x = \left(\matrix  0&1\\
1&0\endmatrix \right) (\bar v_x)^{-1} \left(\matrix  0&1\\
1&0\endmatrix \right)$, it follows that $Q(x)$ is of the form
$\left(\matrix  0&q(x)\\ \ovl{q(x)}&0\endmatrix \right)$ for some
potential $q(x)$. Set $\cl I(r)\equiv q$. The following results
show that $\cl I$ maps $H^{1,1}_1$ into $H^{1,1}(\rb_-)$.

\proclaim{Lemma 3.14}
For $x\le 0$,
$$
\|C^+(I-v_{x-})\|_{L^2}, \|C^-(v_{x+}-I)\|_{L^2} \le
\frac1{(1+x^2)^{1/2}}\|r\|_{H^{1,0}}.
$$
\endproclaim

\demo{Proof} We will only estimate $u \equiv C^-(v_{x+}-I)$ in
which $u_{21}$ is the only nonzero entry; $C^+(I-v_{x-})$ may be
estimated in a similar manner. By Fourier transform
$$
r(z) = \frac1{\sqrt{2\pi}} \int^{+\infty}_{-\infty} e^{i\xi z}\hat r(\xi)d\xi,
$$
and hence
$$
u_{21}(x,\langle z\rangle) = C^- e^{-ix\langle z\rangle} (-\bar r) = \frac{-1}{\sqrt{2\pi}} \int^x_{-\infty} e^{i(\xi-x)\langle z\rangle} \widehat{(-\bar r)}(\xi) d\xi,\tag3.15$$
where $f(\langle z\rangle)$ denotes the function $z\mapsto f(z)$. Thus we have
$$\align
\|u_{21}(x,\cdot)\|_{L^2(dz)} &= \left(\int^x_{-\infty} |\hat{\bar r}(\xi)|^2 d\xi\right)^{\frac12}\tag3.16\\
&\le (1+x^2)^{-\frac12} \left(\int^x_{-\infty} (1+\xi^2) |\hat{\bar r}(\xi)|^2d\xi\right)^{1/2}\\
&\le (1+x^2)^{-\frac12} \|r\|_{H^{1,0}}.
\endalign$$
\enddemo

\proclaim{Theorem 3.17} $\cl I$ maps $H^{1,0}_1$ into $H^{0,1}$.
\endproclaim

\demo{Proof}
As above, we only consider $x\le 0$. Write
$$
\int_{\rb} ((1-C_{v_{x\pm}})^{-1} I)(v_{x+}-v_{x-})dz \equiv
\int_1 + \int_2 + \int_3,\tag3.18$$ where
$$\align
\int_1 &= \int (v_{x+}-v_{x-}),\\
\int_2 &= \int(C_{v_{x\pm}}I)(v_{x+}-v_{x-})\\
\int_3 &= \int (C_{v_{x\pm}}(1-C_{v_{x\pm}})^{-1}C_{v_{x\pm}}I)(v_{x+} - v_{x-})\\
&= \int (C_{v_{x\pm}}(\mu-I))(v_{x+}-v_{x-}).
\endalign$$
We remark that for calculating $q$, the estimate of $\int_2$ is not needed because it is diagonal and $\text{ad } \sigma\int_2 =0$. But the estimate is useful for other problems. Clearly, $\int(v_{x+}-v_{x-}) \in H^{0,1}(dx)$ by the Fourier transform. Using the triangularity of $v_\pm$, the fact that $C^+-C^- = 1$, Cauchy's theorem, and Lemma~3.14, we have for $x\le 0$,
$$\align
\left|\int_2\right| &= \left|\int [(C^+(I-v_{x-}))(v_{x+}-I) + (C^-(v_{x+}-I))(I-v_{x-})]\right|\tag3.19\\
&= \left|\int [(C^+(I-v_{x-})) (C^-(I-v_{x+})) + (C^-(v_{x+}-I)) (C^+(I-v_{x-}))]\right|\\
&\le c(1+x^2)^{-1}, \text{ for some } c>0.
\endalign$$
Finally, by Lemma 3.14 and the uniform $L^2$-boundedness of
$(1-C_{v_{x\pm}})^{-1}$,
$$
\|\mu-I\|_{L^2} = \|(1-C_{v_{x\pm}})^{-1} C_{v_{x\pi}}I\|_{L^2} \le c(1+x^2)^{-1/2} \text{ for some } c>0,\tag3.20$$
we have
$$\align
\left|\int_3\right| &= \left|\int[(C^+(\mu-I)(I-v_{x-})) (v_{x+}-I)
+ (C^-(\mu-I)(v_{x+}-I))(I-v_{x-})]\right|\tag3.21\\
&= \left|\int[(C^+(\mu-I)(I-v_{x-}))C^-(v_{x+}-1) + (C^-(\mu-I)(v_{x+}-I))C^+(I-v_{x-})]\right|\\
&\le c(1+x^2)^{-1}, \text{ some } c>0.
\endalign$$
\enddemo

\proclaim{Theorem 3.22}
If $r\in H^{1,1}_1$, then $q\in H^{1,0}$.
\endproclaim

\demo{Proof}
Using the relation $\partial_x\mu = (iz \text{ ad } \sigma + Q)\mu$ and the fact that $\text{ad } \sigma$ is a derivation, we have
$$
\partial_x(\mu e^{ixz \text{ ad } \sigma}(v_+-v_-)) = (iz \text{ ad } \sigma + Q)(\mu e^{ixz \text{ ad } \sigma}(v_+-v_-)).
$$
Thus
$$
Q' = \frac{\text{ad } \sigma}{2\pi} Q \int \mu(x_{x+} - v_{x-}) +
\frac{\text{ad } \sigma}{2\pi} \int i \text{ ad } \sigma \
\mu(zv_{x+}-zv_{x-})).\tag3.23$$ The first term still gives rise
to an $H^{0,1}$ function. On the other hand, the estimate (3.20)
shows the second term in (3.23) gives rise to an $L^2$
function.$\qquad\square$
\enddemo

Now consider the normalized RHP $(\rb,\check{v}_x =
\delta^{\sigma_3}_-v_x \delta^{-\sigma}_+)$ for $\check{m} =
m\delta^{-\sigma_3}$ as above. By the properties of $\delta =
\delta_{z_0=+\infty}$, one sees that $\check{v}_x$ has the form
$\left(\matrix  1&\check{r}e^{izx}\\ -\bar{\check{r}} e^{-izx}&1 -
|\tilde r(z)|^2\endmatrix \right)$, where $\check{r} =
r\delta_+\delta_-$, $\|\check{r}\|_{L^\infty} = \|r\|_{L^\infty} <
1$. There is a completely parallel theory for this RHP, with the
difference that the problem is now adapted for studying the decay
at $x\to +\infty$. In particular one finds that $\check{m}_\pm$
solve the differential equation $\frac{d}{dx} \check{m}_\pm =
iz(\sigma, \check{m}_\pm) + \check{Q}\check{m}_\pm$, where
$\check{Q}(x) = \left(\matrix  0&\check{q}(x)\\
\ovl{\check{q}(x)}&0\endmatrix \right) \in H^{1,1}(\rb_+)$. But as
$\delta = \delta(z)$ is independent of $x$, it follows that the
differential equations for $m_\pm$ and $\check{m}_\pm$ can only be
compatible if $Q=\check{Q}$ i.e.\ $q(x) = \check{q}(x)$:\ we
conclude in particular that $q\in H^{1,1}(\rb)$. (More precisely,
we have shown that $\check{q}$ extends $q(x)$ defined on $\rb_-$
to an $H^{1,1}$ function on all of $\rb$.)

Let $m(x,z) = I+(C(\mu(w^+_x  + w^-_x)))(z)$ be the extension of
$m_\pm-I$ off $\Sigma=\rb$. Then by analytic continuation, for
each $z\in \cb\backslash \rb$, $m = m(x,z)$ solves the
differential equation $\frac{dm}{dx} = iz[\sigma,m] + Qm$ and by
the Riemann--Lebesgue lemma and (3.20), we conclude that

$\ds\bullet\hskip1in m(x,z)\to I\quad \text{as}\quad x\to -\infty.$\medskip

\n Similarly $\check{m}(x,z)\to I$ as $x\to +\infty$ and hence

 $\ds\bullet\hskip1in m(x,z) = \check{m}(x,z) \delta(z)^{\sigma_3} \quad \text{is bounded as}\quad  x\to +\infty.$\medskip

\n This shows that $m(x,z)$ is the (unique) \cite{BC}-type
solution for $Q = \left(\matrix  0&q(x)\\ \bar q(x)&0\endmatrix
\right)$ for all $z\in \cb \backslash \rb$. As $m_+ = m_- v_x =
m_- \left(\matrix  1-|r|^2&re^{izx}\\ -\bar re^{-izx}&1\endmatrix
\right)$, this means by definition that $\cl R(q) = r$ i.e.\ $\cl
R(\cl I(r)) = r$. On the other hand $\cl R$ is one-to-one, for if
$r=\cl R(q) = \cl R(q^\#)$ then we have two solutions $m$ and
$m^\#$ of the normalized RHP $(\rb, v_x)$, and hence $\Delta m_+ =
\Delta m_- v_x$ on $\rb$, $\Delta m_\pm \in \partial C(L^2)$,
where $\Delta m = m-m^\#$. But again by Proposition 2.14, if
$\Delta m\ne 0$, this implies that $\text{Ker}(1-C_{w_x})\ne
\{0\}$, which is a contradiction. Hence $m(x,z) = m^\#(x,z)$ for
all $x$ and so $q(x) = q^\#(x)$. Thus $\cl R$ is a bijection with
inverse $\cl R^{-1}=\cl I$. This completes the proof that $\cl R$
is bijective from $H^{1,1}$ onto $H^{1,1}_1$.

Finally, as noted above, the solution $m_\pm$ of $(\Sigma,v_x)$
has the form $m_\pm = I+C^\pm(\mu(w^+_x + w^-_x))$. Let $m(x,z) =
I+C(\mu(w^+_x+w^-_x))(z)$ again be the extension of $m_\pm$ off
$\rb$. As $r\in H^{1,1}_1$, we see that
$$
m(x,z) = I + \frac{m_1(x)}z + o\left(\frac1z\right) \tag3.24$$ as
$z\to\infty$ in any proper subsector of $\cb^+$ or $\cb^-$, where
the {\it residue\/} $m_1(x)$ of $m(x,z)$ is given by
$$
m_1(x) = -\frac1{2\pi i} \int_{\rb}\mu(w^+_x + w^-_x)ds.
\tag3.25$$ We conclude in particular that
$$
Q(x) = -i \text{ ad } \sigma(m_1(x))\tag3.26$$
which implies
$$
q(x) = -i(m_1(x))_{12}\tag3.27$$ as in (1.6) above.

\demo{Remark 3.28} The above calculations also show that $\cl R$
is a bijection from $H^{0,1}$ onto $H^{1,0}_1$.

\enddemo

 \head 4. Smoothing estimates and the proof of
Theorem 1.10
\endhead

Throughout this section we always assume that $r\in H^{1,1}_1$,
which corresponds to initial data $q_0 = q(t=0) = \cl R^{-1}(r)$
in $H^{1,1}$ for NLS, by the results of \S 3. The (unique, weak)
solution of NLS in $H^{1,1}$ with initial data $q_0$ is given by
(1.8), $q(t) = \cl R^{-1}(e^{-i\lozenge^2t}r(\lozenge))$. In terms
of the normalized RHP $(\Sigma  = \rb,v_\theta)$,
$$
m_+ = m_-v_\theta \quad \text{on}\quad \rb,\qquad m_\pm \in I +
\partial C(L^2),\tag4.1$$ with $v_\theta = \left(\matrix
1-|r|^2&re^{i\theta}\\ -\bar re^{-i\theta}&1\endmatrix \right)=
(I-w^-_\theta)^{-1} (I+w^+_\theta)$ and
$$
w_\theta = (w^-_\theta, w^+_\theta) = \left(\left(\matrix  0&re^{i\theta}\\ 0&0\endmatrix \right), \left(\matrix  0&0\\ -\bar re^{-i\theta}&0\endmatrix \right)\right),
\tag4.2$$
the solution $q(t)$ is given by (3.13)
$$
Q(t) = \left(\matrix  0&q(t)\\ \ovl{q(t)}&0\endmatrix \right) =
\frac{\text{ad } \sigma}{2\pi} \int_{\rb} \mu(w^+_\theta +
w^-_\theta)dz\tag4.3$$ where $\mu\in I+L^2$ is the unique solution
of $(1-C_{w_\theta})\mu = I$.

As noted before, the steepest descent method proceeds by taking
advantage of the signature table for $\text{Re } i\theta$ (see
Figure~2.53). The key step (cf.\ \cite{DIZ}, \cite{DZ2}) is to
separate the factor $e^{i\theta}$ algebraically from the factor
$e^{-i\theta}$. For $z>z_0$, we use the upper/lower factorization
$$
v_\theta = \left(\matrix  1&re^{i\theta}\\ 0&1\endmatrix \right) \left(\matrix  1&0\\ -\bar re^{-i\theta}&1\endmatrix \right)
$$
and for $z<z_0$, we use the lower/diagonal/upper factorization
$$
v_\theta = \left(\matrix  1&0\\ -\frac{\bar re^{-i\theta}}{1-|r|^2}&1\endmatrix \right) \left(\matrix  1-|r|^2&0\\ 0&\frac1{1-|r|^2}\endmatrix \right) \left(\matrix  1&\frac{re^{i\theta}}{1-|r|^2}\\ 0&1\endmatrix \right).
$$
The diagonal factors $(1-|r|^2)^{\pm1}$ can be removed by
conjugating $v_\theta$ by $\delta^{\sigma_3}_\pm$, where
$\delta_\pm = e^{C^\pm_{\rb_- +z_0} \log 1-|r|^2}$ as in (2.39).
We obtain for $z>z_0$,
$$
\check{v}_\theta \equiv \delta^{\sigma_3}v_\theta \delta^{-\sigma_3} = \left(\matrix  1&r\delta^2 e^{i\theta}\\ 0&1\endmatrix \right) \left(\matrix  1&0\\ -\bar r\delta^{-2}e^{i\theta}&1\endmatrix \right)
$$
and for $z<z_0$
$$
\check{v}_\theta \equiv \delta^{\sigma_3}_- v_\theta \delta^{-\sigma_3}_+ = \left(\matrix  1&0\\ -\frac{\bar r\delta^{-2}_-e^{-i\theta}}{1-|r|^2}&1\endmatrix \right) \left(\matrix  1&\frac{r\delta^2_+e^{i\theta}}{1-|r|^2}\\ 0&1\endmatrix \right)
$$
(cf.\ $(\rb, \check{v}_x)$ in \S 3). Note that the multipliers
$e^{i\theta}$ and $e^{-i\theta}$ have now been separated in the
sense that they lie in different factors of $\check{v}_\theta$. If
$m_\pm$ is the solution of the normalized RHP $(\rb,v_\theta)$,
then $\check{m}_\pm \equiv m_\pm \delta^{-\sigma_3}_\pm$ clearly
solves the normalized RHP $(\rb,\check{v}_\theta)$. Rewriting the
jump relation for $\check{v}_\theta$ for $z>z_0$ in the form
$\check{m}_+ \left(\matrix  1&0\\ -\bar
r\delta^2e^{-i\theta}&1\endmatrix \right)^{-1} = \check{m}_-
\left(\matrix  1&r\delta^2e^{i\theta}\\ 0&1\endmatrix \right)$,
and neglecting, for the moment,  analyticity issues for the
coefficients $r$, $-\bar r$, we see that $\check{m}_1\left(\matrix
1&0\\ -\bar r\delta^2 e^{-i\theta}&1\endmatrix \right)^{-1}$ can
be continued to a sector above $\rb_+ + z_0$ and $\check{m}_-
\left(\matrix 1&r\delta^2e^{i\theta}\\ 0&1\endmatrix \right)$ can
be continued to the sector below $\rb_++z_0$. The same is true for
the appropriate factors on $\rb_-+z_0$. We then obtain a RHP on a
cross $(z_0 + e^{i\pi/4}\rb) \cup (z_0 + e^{-i\pi/4}\rb)$, say,
and things are so arranged so that all the factors $e^{\pm
i\theta}$ are now exponentially decreasing. As $t\to\infty$, the
RHP problem then localizes at $z_0$.

The main analytical task in the method is to handle the lack of
analyticity of the coefficients $r$, $-\bar r$ etc. In \cite{DIZ},
\cite{DZ2}, this is done by approximating these coefficients by
rational functions to high enough order at $z=z_0$, and this
requires a high order of smoothness and decay for $r$. In this
paper we show that a suitable approximation can still be done when
$r$ is just in $H^{1,1}_1$, but now we must utilize cancellations
from oscillations, and not just absolute type estimates. As
mentioned in \S 1, Theorem~1.10 with $r= \cl R(q_0) \in H^{1,1}_1$
is precisely what we need for the perturbation theory of NLS.

\bigskip\bigskip
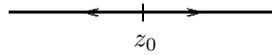
\font\thinlinefont=cmr5 \centerline{\beginpicture
\setcoordinatesystem units <.70000cm,.70000cm>
\linethickness=1pt
\setshadesymbol ({\thinlinefont .}) \setlinear
%
%
\linethickness= 0.500pt \setplotsymbol ({\thinlinefont .})
%
%
\linethickness= 0.500pt \setplotsymbol ({\thinlinefont .})
%
%
\linethickness= 0.500pt \setplotsymbol ({\thinlinefont .})
%
%
%
%
%
\linethickness= 0.500pt \setplotsymbol ({\thinlinefont .})
%
%
%
%
%
\linethickness= 0.500pt \setplotsymbol ({\thinlinefont .})
%
%
%
%
%
\linethickness= 0.500pt \setplotsymbol ({\thinlinefont .})
%
%
%
%
%
\linethickness= 0.500pt \setplotsymbol ({\thinlinefont .})
\putrule from  2.540 21.590 to  7.620 21.590
%
%
\linethickness= 0.500pt \setplotsymbol ({\thinlinefont .})
\putrule from  5.080 21.749 to  5.080 21.431
%
%
\linethickness= 0.500pt \setplotsymbol ({\thinlinefont .})
\putrule from  6.032 21.590 to  6.191 21.590
%
%
\plot  5.937 21.526  6.191 21.590  5.937 21.654 /
%
%
%
\linethickness= 0.500pt \setplotsymbol ({\thinlinefont .})
\putrule from  4.128 21.590 to  3.969 21.590
%
%
\plot  4.223 21.654  3.969 21.590  4.223 21.526 /
%
%
%
\put{$z_0$} [lB] at  4.9 20.955
%
%
\endpicture}

\bigskip
\centerline{Figure 4.4 \ $\tilde \rb_{z_0}$}

For reasons that will become clear further on we reverse the
orientation on $\rb_- +z_0$ to obtain a contour
$\widetilde{\rb}_{z_0} = e^{i\pi}(\rb_++z_0) \cup (\rb_++z_0)$
with associated jump matrix $\tilde v_\theta = \check{v}_\theta$
for $z>z_0$ and $\tilde v_\theta = \check{v}^{-1}_\theta$ for
$z<z_0$ (see discussion preceding Proposition~2.20). Observe that
if $\check{m}_\pm = I + C^\pm_{\rb}  \check{h} \in I + \partial
C(L^2)$ is the solution of the normalized RHP $(\rb,
\check{v}_\theta)$, then $\tilde m_\pm = I +
C^\pm_{\widetilde{\rb}_{z_0}}\tilde h \in I + \partial C(L^2)$ is
the solution of the normalized RHP $(\widetilde{\rb}_{z_0}, \tilde
v_\theta)$ if $\tilde h \equiv \tilde h$ for $z>z_0$ and $\tilde
h\equiv -\check{h}$ for $z<z_0$, and vice versa. Moreover, the
extension $\tilde m(z)$ of $\tilde m_\pm$ of
$\widetilde{\rb}_{z_0}$ is the same as the extension
$\check{m}(z)$ of $\check{m}_\pm$ off $\rb$, and is clearly given
by $m(z)\delta(z)^{-\sigma_3}$, $z\in \cb\backslash \rb$. If $m_1
= m_1(x,t)$ and $\delta_1$ are the residues (cf.\ (3.24)) of $m(z)
= m(x,t,z)$ and $\delta(z)^{-\sigma_3}$ respectively,
$$\align
m(x,t,z) &= I+\frac{m_1(x,t)}z + o\left(\frac1z\right),\\
\delta(z)^{-\sigma_3} &= I + \frac{\delta_1}z +
o\left(\frac1z\right) = I + \left(\frac1{2\pi i} \int_{\rb_- +
z_0} \log (1 - |r|^2)\right) \frac{\sigma_3}z +
o\left(\frac1z\right),
\endalign$$
then we see that the residue $\tilde m_1(x,t)$ of $\tilde
m(x,t,z)$ is given by $\tilde m_1(x,t) = m_1(x,t)  + \delta_1$.
But then by (3.27), as $\delta_1$ is diagonal,
$$
q(x,t) = -i(m_1(x,t))_{12} = -i(\tilde m_1(x,t))_{12}.
\tag4.5$$

From the form of $\check{v}_\theta$, we see that $\tilde v_\theta = (I-\widetilde w^-_\theta)^{-1}(I + \widetilde w^+_\theta)$ where
$$
\widetilde w_\theta = (\widetilde w^-_\theta, \widetilde w^+_\theta) = \left\{\matrix
\left(\left(\matrix  0&r\delta^2e^{i\theta}\\ 0&0\endmatrix\right),
\left(\matrix  0&0\\ -\bar r\delta^{-2}e^{-i\theta}&0\endmatrix \right)\right), \qquad z>z_0,\hfill\\
\left(\left(\matrix  0&-r\delta_+\delta_-e^{i\theta}\\
0&0\endmatrix \right), \left(\matrix  0&0\\  \bar r\delta^{-1}_+
\delta^{-1}_- e^{-i\theta}&0\endmatrix \right)\right),\qquad
z<z_0,
\endmatrix\right.\tag4.6$$
which can also be written as
$$
\widetilde w_\theta = \left(\left(\matrix
0&\frac{-r\tilde\delta^2_- e^{i\theta}}{1-|r|^2}\\ 0&0\endmatrix
\right), \left(\matrix  0&0\\ \frac{\bar r\tilde \delta^{-2}_+
e^{-i\theta}}{1-|r|^2}&0 \endmatrix \right)\right),\qquad
z<z_0,\tag4.7$$ where $\tilde\delta_\pm(z)$ denotes the boundary
values of $\delta(z)$ on $\widetilde{\rb}_{z_0}$. Thus
$\tilde\delta_\pm(z) = \delta_\pm(z)$ for $z>z_0$ and
$\tilde\delta_\pm(z) = \delta_\mp(z)$ for $z<z_0$. Observe that if
we reverse the orientation as above for the normalized RHP
$(\rb,\hat v_\theta)$ (see (2.48) et seq), we obtain a normalized
RHP $(\widetilde{\rb}_{z_0}, \bold{v}_\theta =
(I-\bold{w}^-_\theta)^{-1}(I + \bold{w}^+_\theta))$ where
$\bold{w}_\theta = (\bold{w}^-_\theta, \bold{w}^+_\theta)$ is the
same as $\widetilde w_\theta$ except $r$ is replaced by $[r]$,
etc.

Our goal is to prove that the solution of $(\widetilde{\rb}_{z_0},
\tilde v_\theta)$ is close to the solution of
$(\widetilde{\rb}_{z_0}, \bold{v}_\theta)$ as $t\to\infty$. For
later purposes, note that what the reversal of orientation on
$\rb_-+z_0$ achieves, is that the triangularity of $\widetilde
w^+_\theta$ is the same for $z>z_0$ as for $z<z_0$, and similarly
for $\widetilde w^-_\theta$. Also note that by Proposition~2.20,
$C_{\tilde w_\theta} = C_{\bold{w}_\theta}$, and hence by
Proposition~2.49, $(1 - C_{\bold{w}_\theta})^{-1}$ exists for
large $t$, say $t\ge t_0$, and
$$
\|(1-C_{\bold{w}_\theta})^{-1}\|_{L^2(\tilde{\rb}_{z_0})} \le
c\tag4.8$$ for all $x\in \rb$ and for all $t\ge t_0$.

Also, by (2.30), $\|(1-C_{w_\theta})^{-1}\|_{L^2(\rb)} \le
\frac1{1-\|r\|_{L^\infty}}$, and again using Proposition~(2.20) it
is easy to see that $(1-C_{\tilde w_\theta})^{-1}$ exists in
$L^2(\widetilde{\rb}_{z_0})$ and we have the bound
$$
\|(1-C_{\tilde w_\theta})^{-1}\|_{L^2(\tilde{\rb}_{z_0})} \le
c.\tag4.9$$

First we need some technical lemmas. We use $\rho<1,\lambda$ and $\eta$ to denote $L^\infty$, $H^{1,0}$ and $H^{1,1}$ bounds for $r$ respectively. Thus
$$
\|r\|_{L^\infty} \le \rho,\quad \|r\|_{H^{1,0}} \le \lambda,\quad \|r\|_{H^{1,1}} \le \eta.
$$
Extend $\widetilde{\rb}_{z_0}$ to a contour $\Gamma_{z_0} =
\widetilde{\rb}_{z_0} \cup (z_0 + e^{i\pi/2} \rb_-) \cup (z_0 +
e^{-i\frac\pi2}\rb_-)$.

\bigskip\bigskip
\font\thinlinefont=cmr5 \centerline{\beginpicture
\setcoordinatesystem units <.70000cm,.70000cm>
\linethickness=1pt
\setshadesymbol ({\thinlinefont .}) \setlinear
%
%
\linethickness= 0.500pt \setplotsymbol ({\thinlinefont .})
\putrule from 12.700 21.590 to 17.780 21.590
%
%
\linethickness= 0.500pt \setplotsymbol ({\thinlinefont .})
\putrule from 15.240 24.130 to 15.240 19.050
%
%
\linethickness= 0.500pt \setplotsymbol ({\thinlinefont .})
\putrule from 16.034 21.590 to 16.192 21.590
%
%
\plot 15.939 21.526 16.192 21.590 15.939 21.654 /
%
%
%
\linethickness= 0.500pt \setplotsymbol ({\thinlinefont .})
\putrule from 14.129 21.590 to 13.970 21.590
%
%
\plot 14.224 21.654 13.970 21.590 14.224 21.526 /
%
%
%
\linethickness= 0.500pt \setplotsymbol ({\thinlinefont .})
\putrule from 15.240 23.178 to 15.240 22.860
%
%
\plot 15.176 23.114 15.240 22.860 15.304 23.114 /
%
%
%
\linethickness= 0.500pt \setplotsymbol ({\thinlinefont .})
\putrule from 15.240 20.003 to 15.240 20.161
%
%
\plot 15.304 19.907 15.240 20.161 15.176 19.907 /
%
%
%
\linethickness= 0.500pt \setplotsymbol ({\thinlinefont .})
%
%
\linethickness= 0.500pt \setplotsymbol ({\thinlinefont .})
%
%
\linethickness= 0.500pt \setplotsymbol ({\thinlinefont .})
%
%
%
%
%
\linethickness= 0.500pt \setplotsymbol ({\thinlinefont .})
%
%
%
%
%
%
%
\put{$z_0$} [lB] at 14.6 21.114
\endpicture}

\bigskip
\centerline{Figure 4.10\ \ \ \ $\Gamma_{z_0}$}\bigskip

\n As a {\it complete\/} contour (see e.g. \cite{Z1}), $\Gamma_{z_0}$ has the important property
$$
C^+_{\Gamma_{z_0}} C^-_{\Gamma_{z_0}} = C^-_{\Gamma_{z_0}}
C^+_{\Gamma_{z_0}} = 0.\tag4.11$$ In the following Lemmas we will
assume for convenience that $x=0$. Thus $z_0=0$, $\delta =
\delta_{z_0=0}$, $\Delta = \Delta_{z_0=0}$, $\widetilde{\rb}
\equiv \widetilde{\rb}_{z_0=0}$, $\Gamma \equiv
\widetilde\Gamma_{z_0=0}$, and $\theta = -tz^2$. Observe that the
bounds in the Lemmas depend only on $\rho$ and $\lambda$, but not
on $\eta$.

\proclaim{Lemma 4.12} For $z\in \rb \backslash 0$,
$$
|\Delta'(z)|\le \text{\rm I} + \text{\rm II},\tag4.13$$ where
$$\align
\|\text{\rm I}\|_{L^2} &\le \frac{c\rho\lambda}{1-\rho}
\tag4.14\\
\text{\rm II} &\le \frac{c\rho^2}{1-\rho} \frac1{|z|}.
\tag4.15
\endalign$$
\endproclaim

\demo{Proof}
We have
$$\align
\Delta'(z) &= -\Delta(z) \frac{d}{dz} H((\log(1-|r|^2)) \chi_{\rb_-})\\
&= \Delta(z)H \left(\frac{|r|^{2'}}{1-|r|^2} \chi_{\rb_-}\right) -
\frac{i\Delta}\pi \frac{\log(1-|r(0)|^2)}{z}.
\endalign$$
The result now follows from the $L^2$ mapping properties of $H$,
the identity $|\Delta| = 1$, and the elementary bound
$|\log(1-|r(0)|^2)| \le \frac{|r(0)|^2}{1-|r(0)|^2}$.
$\qquad\square$\enddemo

\proclaim{Lemma 4.16}
Suppose $f\in H^{1,0}$. Then for all $t>1$,
$$
\left|\int_{\rb} f\Delta^{\pm1} e^{\mp itz^2}\ dz\right| \le
\frac{c}{t^{1/2}} \frac{1+\lambda}{(1-\rho)}
\|f\|_{H^{1,0}}.\tag4.17$$ In addition, if $f(0) = 0$, then for
all $t>1$,
$$
\left|\int_{\rb} f\Delta^{\pm 1}e^{\mp itz^2} \ dz\right| \le
\frac{c}{t^{3/4}} \frac{1+\lambda}{(1-\rho)}
\|f\|_{H^{1,0}}.\tag4.18$$
\endproclaim

\demo{Proof}
We only consider the case $\Delta = \Delta^{+1}$ above. The other case is similar. Decompose the  integral as follows:
$$
\int^\infty_{-\infty} f\Delta e^{-itz^2}\ dz = \int_{|z|< \frac1{\sqrt t}} f\Delta e^{-itz^2}\ dz + \int_{|z|> \frac1{\sqrt t}} f\Delta e^{-itz^2}\ dz \equiv \text{I} + \text{II}.
$$
Changing variables, we obtain
$$
|\text{I}| = \frac1{\sqrt t} \left|\int^1_{-1} f\left(\frac{z}{\sqrt t}\right) \Delta\left(\frac{z}{\sqrt t}\right) e^{-iz^2} \ dz\right| \le \frac1{\sqrt t} \|f\|_{L^\infty} \le \frac{c}{\sqrt t} \|f\|_{H^{1,0}}.
$$
We consider $z<-t^{-1/2}$. The case $z>t^{1/2}$ is similar. Integration by parts leads to
$$\align
\int_{z<-t^{-1/2}} f\Delta e^{-itz^2}\ dz &= \frac{e^{-i}}{2i\sqrt t} f \left(-\frac1{\sqrt t}\right) \Delta \left(- \frac1{\sqrt t}\right) + \frac1{2it} \int^{-t^{-1/2}}_{-\infty} e^{-itz^2} \left(\frac{f'\Delta}z + \frac{f\Delta'}z - \frac{f\Delta}z\right)dz\\
&\equiv \text{I}' + \text{II}'' + \text{III}' + \text{IV}'.
\endalign$$
Clearly
$$
|\text{I}'| \le \frac{c\|f\|_{H^{1,0}}}{\sqrt t}, \qquad |\text{II}'| \le \frac{c}{t^{3/4}} \|f\|_{H^{1,0}}, \quad \text{and}\quad |\text{IV}'| \le \frac{c\|f\|_{H^{1,0}}}{t^{1/2}}.
$$
Finally using Lemma (4.12), we obtain
$$
|\text{III}'| \le \frac{\|f\|_{L^\infty}}{2t} \int^{-t^{-1/2}}_{-\infty} \left|\frac{\Delta'}z\right| dz \le \frac{\|f\|_{H^{1,0}}}{2t} \left(\frac{\rho\lambda}{1-\rho} t^{1/4} + \frac{c\rho^2}{1-\rho} t^{1/2}\right)
$$
which now leads directly to (4.17). If $f(0) = 0$, then the same
arguments together with the bound $|f(z)|\le
|z|^{1/2}\|f\|_{H^{1,0}}$ for $|z|\le 1$, say, in I, I$'$, III$'$,
and IV$'$, yield (4.18). $\qquad\square$\enddemo

Let $D_j, j=1,\ldots, 4$, be the $j^{\text{\rm th}}$ quadrant in
$\cb\backslash\Gamma$

$$\vbox{\offinterlineskip \halign{\strut \quad \hfil $#$\hfil\quad &
\vrule\quad \hfil $#$\hfil\quad\cr D_2&D_1\cr \noalign{\hrule}
D_3&D_4\cr}}$$ \centerline{Figure~4.19\ \ \ $\Gamma$}\bigskip

\n In the Lemma below $H^q$ denotes Hardy space. A general
reference for Hardy spaces is, for example, \cite{Dur}.

\proclaim{Lemma 4.20}
Suppose $f\in H^{1,0}$, then for $2\le p < \infty$ and for all $t\ge 0$,
$$
\left\{\matrix
\|C^-_{\rb_+\to\Gamma} \delta^{-2} fe^{it\lozenge^2}\|_{L^p} \le {\dsize\frac{c}{(1+t)^{1/2p}}} \|\delta^{-2}\|_{L^\infty(D_1)} \|f\|_{H^{1,0}} \le {\dsize\frac{c}{(1+t)^{1/2p}}} {\dsize\frac{\|f\|_{H^{1,0}}}{1-\rho}},\hfill\\
\|C^-_{e^{i\pi} \rb_+\to \Gamma} \tilde \delta^{-2}_+ fe^{it\lozenge^2} \|_{L^p} \le {\dsize\frac{c}{(1+t)^{1/2p}}} \|\delta^{-2}\|_{L^\infty(D_3)} \|f\|_{H^{1,0}} \le {\dsize\frac{c}{(1+t)^{1/2p}}} \|f\|_{H^{1,0}},\hfill\\
\|C^+_{\rb_+\to\Gamma} \delta^2 fe^{-it\lozenge^2}\|_{L^p} \le {\dsize\frac{c}{(1+t)^{1/2p}}} \|\delta^2\|_{L^\infty(D_4)} \|f\|_{H^{1,0}} \le {\dsize\frac{c}{(1+t)^{1/2p}}} {\dsize\frac{\|f\|_{H^{1,0}}}{1-\rho}},\hfill\\
\|C^+_{e^{i\pi}\rb_+\to \Gamma} \tilde\delta^2_-
fe^{-it\lozenge^2}\|_{L^p} \le {\dsize\frac{c}{(1+t)^{1/2p}}}
\|\delta^2\|_{L^\infty(D_2)} \|f\|_{H^{1,0}} \le
{\dsize\frac{c}{(1+t)^{1/2p}}} \|f\|_{H^{1,0}}.\hfill
\endmatrix\right.\tag4.21$$
Suppose in addition that $f(0) = 0$ and that $g$ is a function in
the Hardy space $H^q(\cb\backslash \rb)$ for some $2\le q\le
\infty$. Then for all $t\ge 0$,
$$
\left\{\matrix
\|C^-_{\rb_+\to\Gamma} g_+ fe^{it\lozenge^2}\|_{L^2} \le {\dsize\frac{c}{(1+t)^{\frac12-\frac1q}}} \|g\|_{H^q(\cb \backslash \rb)} \|f\|_{H^{1,0}},\hfill\\
\|C^-_{e^{i\pi}\rb_+\to\Gamma} \tilde g_+ fe^{it\lozenge^2}\|_{L^2} \le {\dsize\frac{c}{(1+t)^{\frac12-\frac1q}}} \|g\|_{H^q(\cb\backslash\rb)} \|f\|_{H^{1,0}},\hfill\\
\|C^+_{\rb_+\to\Gamma} g_-fe^{-it\lozenge^2}\|_{L^2} \le
{\dsize\frac{c}{(1+t)^{\frac12-\frac1q}}}\|q\|_{H^2(\cb
\backslash\rb)} \|f\|_{H^{0,1}},\hfill\\
\|C^+_{e^{i\pi}\rb_+\to\Gamma} \tilde g_- fe^{-it\lozenge^2}
\|_{L^2} \le {\dsize\frac{c}{(1+t)^{\frac12-\frac1q}}}
\|g\|_{H^q(\cb\backslash\rb)}
\|f\|_{H^{1,0}},\hfill\endmatrix\right.\tag4.22$$ where $g_\pm$
are the boundary values of $g$ on $\rb$ and $\tilde g_\pm = g\mp$
on $e^{i\pi}\rb_+$.
\endproclaim

\demo{Proof}
Consider the first inequality in (4.21). The other cases in (4.21) are similar. By Fourier theory, $f(z) = {\ds \frac1{\sqrt{2\pi}} \int} e^{-iyz} \check{f}(y)dy$.
We have for any $\epsilon>0$,
$$
C^-_{\rb_+\to\Gamma} \delta^{-2} e^{-\epsilon\lozenge}
fe^{it\lozenge^2} = \frac1{\sqrt{2\pi}} \int \check{f}(y)
e^{-i\frac{y^2}{4t}} F_1\ dy + \frac1{\sqrt{2\pi}} \int
\check{f}(y) e^{-i\frac{y^2}{4t}} F_2\ dy
$$
where
$$\align
F_1(y) &= C^-_{\rb_+\to \Gamma} (\delta^{-2} e^{-\epsilon
\lozenge} \chi_{(0,a)} e^{it(\lozenge-y/2t)^2})
\tag4.23\\
F_2(y) &= C^-_{\rb_+\to \Gamma} (\delta^{-2} e^{-\epsilon\lozenge}
\chi_{(a,\infty)} e^{it(\lozenge-y/2t)^2})\tag4.24
\endalign$$
and $a = \max(0,y/2t)$. (The factor $e^{-\epsilon z}$ is included
just to ensure that $F_2(y)$ exists in $L^p$.) Clearly $F_1(y)$ is
supported on $\rb_+$. Assume first that $p>2$. Then for $y>0$,
$$
\|F_1\|_{L^p(\Gamma)} \le c\|\delta^{-2}\|_{L^\infty(\rb_+)}
\|\chi_{(0,a)}\|_{L^p} \le c\|\delta^{-2}\|_{L^\infty(\rb_+)}
\frac{|y|^{1/p}}{t^{1/p}}
$$
and hence
$$
\left\|\frac{1}{\sqrt{2\pi}} \int\ dy \check{f}(y) e^{-iy^2/4t}
F_1\right\|_{L^p(\Gamma)} \le
\frac{c\|\delta^{-2}\|_{L^\infty(\rb_+)}}{t^{1/p}} \int^\infty_0
|\check{f}(y)|y^{\frac1p}\ dy \le
\frac{c\|\delta^{-2}\|_{L^\infty(D_1)}}{t^{1/p}}
\|f\|_{H^{1,0}}.\tag4.25$$ For $p=2$, rewrite the integral as
$$
\frac1{\sqrt{2\pi}} C^-_{\rb_+\to\Gamma}
\left(\int^\infty_{2t\lozenge} \delta^{-2} e^{-\epsilon \lozenge}
\check{f}(y) e^{-iy\lozenge} e^{it\lozenge^2} \ dy\right).
$$
Using Hardy's inequality \cite{HLP}
$$
\left\|\int^\infty_{2t\lozenge}
|\check{f}(y)|dy\right\|_{L^2(\rb_+)} \le \sqrt{\frac2t}
\|\lozenge \check{f}\|_{L^2(\rb_+)}\tag4.26$$ and hence
$$
\left\|\frac1{\sqrt{2\pi}} \int \check{f} e^{-\epsilon \lozenge} e^{-iy^2/4t} F_1\ dy\right\|_{L^2(\Gamma)} \le \frac{c\|\delta^{-2}\|_{L^\infty(D_1)}}{t^{1/2}} \|f\|_{H^{1,0}}.
$$
For $F_2$, first consider the case when $y<0$, and hence $a=0$. Then for $p\ge 2$,
$$\align
\|C^-_{\rb_+\to\Gamma} (\delta^{-2} e^{-\epsilon \lozenge} e^{it(\lozenge-y/2t)^2})\|_{L^p} &= \|C_{e^{i\pi/4} \rb_+\to\Gamma} (\delta^{-2} e^{-\epsilon\lozenge} e^{it(\lozenge-y/2t)^2})\|_{L^p}, \quad \text{by Cauchy's Theorem,}\\
&\le c\|\delta^{-2}\|_{L^\infty(D_1)} \|e^{it(\lozenge-y/2t)^2})\|_{L^p(e^{i\pi/4}\rb_+)}\\
&\le c\|\delta^{-2}\|_{L^\infty(D_1)}\|e^{it\lozenge^2} e^{-iy\lozenge}\|_{L^p(e^{i\pi/4}\rb_+)}\\
&\le c\|\delta^{-2}\|_{L^\infty(D_1)} \|e^{it\lozenge^2}\|_{L^2(e^{i\pi/4}\rb_+)}, \quad \text{as } y<0\\
&\le \frac{c}{t^{1/2p}} \|\delta^{-2}\|_{L^\infty(D_1)}, \quad \text{by scaling.}
\endalign$$
If $y>0$, then $a=y/2t$ and for $p\ge 2$,
$$\align
\|C^-_{\rb_+\to\Gamma} \delta^{-2} e^{-\epsilon\lozenge} \chi_{(a,\infty)} e^{it(\lozenge-y/2t)^2}\|_{L^p(\Gamma)} &= \|C^-_{(y/2t,\infty)\to\Gamma} \delta^{-2} e^{-\epsilon\lozenge} e^{it(\lozenge-y/2t)^2}\|_{L^p(\Gamma)}\\
&= \|C_{(y/2t+e^{i\pi/4}\rb_+)\to \Gamma} \delta^{-2} e^{-\epsilon\lozenge} e^{it(\lozenge-y/2t)^2}\|_{L^p(\Gamma)}, \text{ again by Cauchy,}\\
&\le c\|\delta^{-2}\|_{L^\infty(D_1)}\|e^{it(\lozenge-y/2t)^2} \|_{L^p(y/2t+e^{i\pi/4}\rb_+)}\\
&\le \frac{c\|\delta^{-2}\|_{L^\infty(D_1)}}{t^{1/2p}}, \text{ again by scaling.}
\endalign$$
Thus for $p\ge 2$,
$$
\left\|\frac1{\sqrt{2\pi}} \int\ dy \check{f}(y) e^{-iy^2/4t} F_2\right\|_{L^p(\Gamma)} \le \frac{c\|\delta^{-2}\|_{L^\infty(D_1)}}{t^{\frac1{2p}}} \|f\|_{H^{1,0}}.
\tag4.27$$

Letting $\epsilon \downarrow 0$ in (4.25)--(4.27), we obtain for $p\ge 2$ and $t\ge 1$,
$$
\|C^-_{\rb_+\to\Gamma} (\delta^{-2}
fe^{it\lozenge^2})\|_{L^p(\Gamma)} < \frac{c}{t^{\frac1{2p}}}
\|\delta^{-2}\|_{L^\infty(D_1)} \|f\|_{H^{1,0}}.
$$
As
$$
\|C^-_{\rb_+\to\Gamma} (\delta^{-2}
fe^{it\lozenge^2})\|_{L^p(\Gamma)} \le
c\|\delta^{-2}\|_{L^\infty(D_1)} \|f\|_{H^{1,0}} \quad \text{for
all} \quad t>0,
$$
we obtain (4.21).

Now we prove the first inequality in (4.22). Again, the remaining inequalities are similar. As before, we have the representation,
$$
C^-_{\rb_+\to\Gamma} g_+ e^{-\epsilon\lozenge} fe^{it\lozenge^2} =
\frac1{\sqrt{2\pi}} \int \check{f}(y) e^{-i\frac{y^2}{4t}} F_1\ dy
+ \frac1{\sqrt{2\pi}} \int \check{f}(y) e^{-i\frac{y^2}{4t}} F_2\
dy
$$
but now as $\int \check{f}(y)dy = \sqrt{2\pi}\ f(0) = 0$,
$$\align
F_1(y) &= C^-_{\rb_+\to\Gamma} (g_+ e^{-\epsilon\lozenge} \chi_{(0,a)} e^{it(\lozenge-y/2t)^2} (1-e^{iy\lozenge}))\\
F_2(y) &= C^-_{\rb_+\to\Gamma} (g_+ e^{-\epsilon \lozenge}
\chi_{(a,\infty)} e^{it(\lozenge-y/2t)^2} (1-e^{iy\lozenge})).
\endalign$$
As before, for $y>0$, the integral involving $F_1$ can be rewritten as
$$
\frac1{\sqrt{2\pi}} C^-_{\rb_+\to\Gamma}
\left(\int^\infty_{2t\lozenge} g_+ e^{-\epsilon\lozenge}
\check{f}(y) (e^{-iy\lozenge} - 1) e^{it\lozenge^2} \ dy\right).
$$
For $q'\ge 2$, using the inequality
$$
\left(\int^\infty_{2tz} |\check{f}(y)|dy\right)^{q'} \le \left(\int^\infty_0 |\check{f}(y)|dy\right)^{q'-2} \left(\int^\infty_{2tz} |\check{f}(y)|dy\right)^2
$$
together with the above Hardy inequality, we obtain for $t>0$
$$
\left\|\int^\infty_{2t\lozenge} |\check{f}(y)|dy
\right\|_{L^{q'}(\rb_+)} \le \frac{c}{t^{1/q'}} \|f\|_{H^{1,0}}.
$$
Hence for $\frac1{q'} + \frac1q = \frac12$,
$$
\left\|\frac1{\sqrt{2\pi}} \int \check{f} e^{-\epsilon \lozenge}
e^{-iy^2/4t} F_1\ dy \right\|_{L^2(\Gamma)} \le
\frac{c\|g_+\|_{L^q(\rb_+)}}{t^{1/q'}} \|f\|_{H^{1,0}} \le
\frac{c\|g\|_{H^q(\cb\backslash \rb)}}{t^{\frac12-\frac1q}}
\|f\|_{H^{1,0}}.
$$
Now for $y<0$, as above
$$
\|F_2(y)\|_{L^2} \le c\|ge^{it(\lozenge-y/2t)^2}
(1-e^{iy\lozenge})\|_{L^2(e^{i\pi/4}\rb_+)} =
c\|ge^{it\lozenge^2}(e^{-iy\lozenge}-1)\|_{L^2(e^{i\pi/4} \rb_+)}.
$$
Hence for $\alpha>0$ and $t>0$,
$$\align
\left\|\frac1{\sqrt{2\pi}} \int^0_{-\infty} \check{f} e^{-\epsilon\lozenge} e^{-iy^2/4t} F_2\ dy\right\|_{L^2(\Gamma)} &\le c \int^0_{-\infty} dy|\check{f}(y)| \left(\int^\infty_0 d\gamma|g(e^{i\pi/4}\gamma)|^2 e^{-2t\gamma^2} |e^{-ie^{i\pi/4}\gamma y} - 1|^2\right)^{\frac12}\\
&\le c \int^0_{-t^\alpha} dy\ |y\check{f}(y)| \left(\int^\infty_0 d\gamma\ \gamma^2|g(e^{i\pi/4}\gamma)|^2 e^{-2t\gamma^2}\right)^{\frac12}\\
&\quad + \int^{-t^\alpha}_{-\infty} dy|\check{f}(y)| \left(\int^\infty_0 d\gamma|g(e^{i\pi/4}\gamma)|^2 e^{-2t\gamma^2}\right)^{\frac12}\\
&\le c(t^{-3/4+\alpha/2} + t^{-\alpha/2-1/4}) \|f\|_{H^{1,0}}
\|g(\lozenge/\sqrt t)\|_{H^q(\cb\backslash \rb)}.
\endalign$$
Taking $\alpha=1/2$, we obtain
$$
\left\|\frac1{\sqrt{2\pi}} \int \check{f} e^{-\epsilon \lozenge}
e^{-iy^2/4t} F_2\ dy\right\|_{L^2(\Gamma)} \le
\frac{c}{t^{1/2-1/(2q)}} \|f\|_{H^{1,0}} \|g\|_{H^q(\cb\backslash
\rb)}.
$$
Finally for $y>0$,
$$
\|F_2(y)\|_{L^2} \le c\|ge^{it(\lozenge-y/2t)^2}
(1-e^{iy\lozenge})\|_{L^2(y/2t+e^{i\pi/4}\rb_+)},
$$
and hence for $\alpha>0$ and $t>0$,
$$\align
&\quad\left\|\frac1{\sqrt{2\pi}} \int^\infty_0 \check{f} e^{-\epsilon\lozenge} e^{-iy^2/4t} F_2\ dy\right\|_{L^2(\Gamma)}\\
& \le c \int^\infty_0 dy|\check{f}(y)| \left(\int^\infty_0 d\gamma|g(y/2t+e^{i\pi/4}\gamma)|^2 e^{-2t\gamma^2} |1-e^{iy(y/2t+e^{i\pi/4}\gamma}|^2\right)^{\frac12}\\
&\le c\int^{t^\alpha}_0 dy|\check{f}(y)| \left(\int^\infty_0 d\gamma|g(y/2t + e^{i\pi/4}\gamma)|^2 e^{-2t\gamma^2} (|e^{-iy^2/2t} -1|^2 + |1-e^{ie^{i\pi/4}\gamma y}|^2)\right)^{\frac12}\\
&\hskip1.5in + c\int^\infty_{t^\alpha} dy|\check{f}(y)| \left(\int^\infty_0 d\gamma|g(y/2t + e^{i\pi/4} \gamma)|^2 e^{-2t\gamma^2}\right)^{\frac12}\\
&\le c(t^{-5/4+3\alpha/2} + t^{-3/4+\alpha/2} +
t^{-\alpha/2-1/4})\|f\|_{H^{1,0}} \|g(\lozenge/\sqrt
t)\|_{H^q(\cb\backslash \rb)}.
\endalign$$
Again, setting $\alpha = 1/2$, we find for $t>0$,
$$
\left\|\frac1{\sqrt{2\pi}} \int^\infty_0 \check{f}
e^{-\epsilon\lozenge} e^{-iy^2/4t} F_2\ dy\right\|_{L^2(\Gamma)}
\le \frac{c}{t^{1/2-1/(2q)}} \|f\|_{H^{1,0}}
\|g\|_{H^q(\cb\backslash \rb)}.
$$
On the other hand, as before,
$$
\|C^-_{\rb_+\to\Gamma}(g_+fe^{it\lozenge^2})\|_{L^2(\Gamma)} \le
c\|g_+\|_{L^q(\rb_+)} \|f\|_{H^{1,0}} \le
c\|g\|_{H^q(\cb\backslash\rb)} ||f\|_{H^{1,0}} \quad \text{for
all}\quad t>0,
$$
and (4.22) follows.$\qquad\square$
\enddemo

\proclaim{Corollary 4.28 to Lemma 4.20}
For any $2\le p < \infty$, and for all $t\ge 0$
$$
\|C^\pm_{\widetilde{\rb}_{z_0}\to \Gamma_{z_0}} \widetilde
w^\pm_\theta\|_{L^p(\Gamma_{z_0})} \le \frac{c}{(1+t)^{1/2p}}
\frac\lambda{(1-\rho)^2}. \tag4.29$$
\endproclaim

\demo{Proof} Translate $z\to z_0 +z$, and then apply Lemma (4.20)
to appropriate choices of f (see $\widetilde w_\theta$ in
(4.6)--(4.7)). The inequality then follows as the $H^{1,0}$ and
$L^\infty$ norm of $r$ is invariant under translation, $r(\cdot)
\to r(\cdot+z_0)$. $\qquad\square$\enddemo

We need the following $L^p$ bound on solutions of RHP's of type
(1.4).

\proclaim{Proposition 4.30} Suppose that r is a continuous
function on $\rb,\lim_{z\to\infty}r(z)=0$ and $\Norm\infty \rb {r}
< \rho <1.$ Then for any $p\ge2$, there exists $t_0=t_0(r,p)$ such
that for $t\ge t_0$ and all $x\in\rb$, $(1-C_{w_\theta})^{-1}$
exists in $L^p(\rb)$ and
$$\|(1-C_{w_\theta})^{-1}\|_{L^p(\rb)\rightarrow L^p(\rb)}\le
c \tag4.31$$ for $t\ge t_0$ and all $x\in\rb$.

\endproclaim

\demo{Remark} In [DZ6] (see also WEBPAGE to [DZ5]) the authors
prove the following stronger, {\it a priori\/} estimate, which is
needed for the perturbation theory in [DZ5]: Suppose $r\in
H_1^{1,0}, \|r\|_{H^{1,0}}\le\lambda,\norm\infty {r} \le \rho<1.$
Then for any $2<p<\infty$, there exists $l_1=l_1(p),l_2=l_2(p)>0$,
and a constant $c=c_p$, such that
$$\|(1-C_{w_\theta})^{-1}\|_{L^p(\rb)\rightarrow L^p(\rb)}\le
c\frac{(1+\lambda)^{l_1}}{(1-\rho)^{l_2}}$$ for all $x,t\in\rb$.
Inequality (4.31, however, is sufficient for the result in this
paper.

In the linear case, the Cauchy operator $C^+$, for example, maps
$L^p\rightarrow L^p, \Norm p \rb {C^+f}\le c_p\Norm p \rb {f}$,
for $1<p<\infty$. In particular, if $\phi=\phi(z)$ is a real
valued function, we see that $\Norm p \rb {C^+fe^{i\theta}}\le
c_p\Norm p \rb {fe^{i\theta}}=c_p\Norm p \rb {f}$, and the bound
is independent of $\phi$. The above bounds on
$(1-C_{w_\theta})^{-1}$ in $L^p$, which are uniform in the
multiplier $e^{i\theta}$, should be viewed as nonlinear versions
of such estimates.

\enddemo

\demo{Proof of Proposition 4.30}

For $\epsilon_1>0$, to be determined below, choose a rational
function $R=R(z)$ such that $$\norm\infty{r-R}\le\epsilon_1$$ and
$$\norm\infty{R}<\rho, \quad \lim_{z\rightarrow\infty}R(z)=0.$$

Associated with $R$ is the model RHP with jump matrix
$v^{\#}=v^{\#}_R$ on $\Sigma^e$ as in (2.61), except $r(z_0)$ is
replaced by $R(z_0)$ etc., and in the definition of $\delta_0$
(see (2.58)) one must use $R(z)$ in place of $r(z)$. Now the
solution $m^{\#}_{\pm}\in I+\dc2$ of the normalized RHP
$(\Sigma^e, v^{\#}=v^{\#}_R)$ can be computed explicity in terms
of parabolic cylinder functions (see [DZ1],[DIZ],[DZ2]) and one
learns that $m_{\pm}^{\#},(m_{\pm}^{\#})^{-1}\in
L^{\infty}(\Sigma^e)$ with bounds which depend {\it only} on
$\rho<1$. Now we know from $\S2$ that the IRHP2$_{L^2(\Sigma^e)},
\quad M_+=M_-v^{\#} + F$, has a unique solution
$M_{\pm}\in\partial C(L^2(\Sigma^e))$, for each $F\in L^2$. But as
$v^{\#}=(m_-^{\#})^{-1}m_+^{\#}$, we have $M_+(m_+^{\#})^{-1}=
M_-(m_-^{\#})^{-1} + F(m_+^{\#})^{-1}$ and hence $M_{\pm}=
(C^{\pm}(F(m_+^{\#})^{-1}))m_{\pm}^{\#}$. But as
$m_{\pm}^{\#},(m_{\pm}^{\#})^{-1}\in L^{\infty}(\Sigma^e)$, it
follows that if $F\in L^2\cap L^p(\Sigma^e)$, then $M_{\pm}\in\dc
p$ for any $F\in L^p$, by density. Moreover, $\norm p {M_{\pm}}\le
c \norm p {F}$, where the constant depends only on $\rho$ and $p$.
Again by Proposition 2.14, this implies $\Norm p
{\Sigma^e}{(1-C_{w^{\#}})^{-1}}\le\frac c{(1-\rho)^2}$ for all
$x,t\in\rb$.

Now let $W_\theta=(W_\theta^-,W_\theta^+)=\left(
\pmtwo0{Re^{i\theta}}00,\pmtwo00{-\bar Re^{-i\theta}}0\right)$ as
in (4.2), but with $r$ replaced by $R$. Then it follows by
arguments similar to the deformation arguments in $\S2$, that the
above bound on $(1-C_{w^{\#}})^{-1}$  implies a similar bound on
$(1-C_{W_\theta})^{-1}$ on $L^p(\rb)$, $\Norm p
{\rb}{(1-C_{W_\theta})^{-1}}\le\frac {c'}{(1-\rho)^2}$, for $t$
sufficiently large, say $t\ge t_0$, and all $x\in \rb$. Again $c'$
depends only on $\rho$ and $p$.

The Proposition follows by making the following choices. First
choose $\epsilon_1>0$ sufficiently small and $R$ as above so that
$\|C_{w_\theta}-C_{W_\theta}\|_{L^p(\rb)\rightarrow
L^p(\rb)}\le\frac12\frac{(1-\rho)^2}{c'}$. Then choose $t_0$
sufficiently large to ensure the a bound $\Norm p
{\rb}{(1-C_{W_\theta})^{-1}}\le\frac {c'}{(1-\rho)^2}$ for $t\ge
t_0$ and all $x\in\rb$. Then use the second resolvent identity to
obtain the bound on $\Norm p {\rb}{(1-C_{w_\theta})^{-1}}$.
$\qquad\square$

\enddemo

\demo{Remark}

The time $t_0$ depends on the full function $r=r(z)$ and not just
on $\norm\infty{r}$. In particular it depends on the location of
the poles of the approximating rational function $R(z)$. If one
assumes that $r\in H_1^{1,0},  \norm\infty{r}\le\rho<1$, as in the
previous Remark, then we may approximate $r$ via the Poisson
formula $R(z)=\int_{\rb}\frac {r(s)}{(s-z)^2+\gamma^2}\frac{\gamma
ds}{\pi}$, leading to bounds $\norm\infty{r-R}\le
c\sqrt{\gamma}\norm2{r'},\|R\|_{H^{1,0}}\le\|r\|_{H^{1,0}}\le\lambda$.
Choosing $\gamma$ sufficiently small, depending only on $\rho$ and
$\lambda$, the reader may check that the proof of the Proposition
now implies that $t_0$ may be chosen to depend only on
$\rho,\lambda$ and $p$, $t_0=t_0(\rho,\lambda,p)$ (for more
details, see [DZ6]).\enddemo

 Reversing the orientation on $\rb$ as above, we conclude that for $r\in H^{1,0}_1$ and $2\le p<\infty$,
$$
\|(1-C_{\widetilde
w_\theta})^{-1}\|_{L^p(\widetilde{\rb}_{z_0})\to
L^p(\widetilde{\rb}_{z_0})} \le c_p \tag4.32$$ where $c_p$ is
uniform for all $x\in\rb$ and all $t\ge t_0$.

\proclaim{Corollary 4.33 to Lemma 4.20}
For any $2\le p < \infty$ and for all $t\ge 0$
$$
\|\tilde \mu-I\|_{L^p} \le \frac{c_p}{(1+t)^{1/2p}}
\frac{\lambda}{(1-\rho)^2}. \tag4.34$$
\endproclaim

\demo{Proof}
The result follows from the formula
$$
\tilde\mu - I = (1-C_{\tilde w_\theta})^{-1} (C_{\tilde w_\theta}I)
$$
together with (4.29) and (4.31).$\qquad\square$
\enddemo

Of course, in order to prove (4.33) all we need are the mapping
properties of $C^\pm_{\widetilde{\rb}_{z_0} \to
\widetilde{\rb}_{z_0}}$, etc. The full estimates in (4.29) for
$C^\pm_{\widetilde{\rb}_{z_0}\to \Gamma_{z_0}}$ are needed below.

\proclaim{Lemma 4.35} For any $2\le p < \infty$ and for $t$
sufficiently large, $\|C^\pm\tilde \mu(\bold{w}_\theta{}^\mp -
\widetilde w^\mp_\theta)\|_{L^2} \le ct^{-\frac12+\frac1{2p}}$.
\endproclaim

\demo{Proof}
By triangularity,
$$\align
C^\pm \widetilde\mu(\bold{w}_\theta{}^\mp -\widetilde w^\mp_\theta) &= C^\pm(\bold{w}_\theta{}^\mp - \widetilde w^\mp_\theta) + C^\pm(C_{\tilde w_\theta}\tilde\mu) (\bold{w}_\theta{}^\mp - \widetilde w^\mp_\theta)\\
&= C^\pm(\bold{w}_\theta{}^\mp - \widetilde w^\mp_\theta) + C^\pm(C^\mp \tilde \mu\widetilde w^\pm_\theta) (\bold{w}_\theta{}^\mp - \widetilde w^\pm_\theta) \equiv \text{I} + \text{II}.
\endalign$$
After translation $z\to z_0+z$, we obtain by (4.22),
$\|\text{I}\|_{L^2}\le ct^{-1/2}$. For $\text{II},\
C^\mp\tilde\mu \widetilde w^\pm_\theta = C^\mp(\tilde
\mu-I)\widetilde w^\pm_\theta + C^\mp \widetilde w^\pm_\theta
\equiv g_1+g_2$. By (4.34), $\|g_1\|_{L^p} \le ct^{-1/2p}$, and by
Lemma~4.20, $\|g_2\|_{L^p} \le ct^{-1/2p}$. It now follows from
(4.22) that $\|\text{II}\|_{L^p} \le ct^{-\frac12+\frac1{2p}}$ and
the lemma is proved.$\qquad\square$
\enddemo

For convenience we abuse the notation by writing $\widetilde w =
\widetilde w^+_\theta + \widetilde w^-_\theta$ and $\bold{w} =
\bold{w}_\theta{}^+ + \bold{w}_\theta{}^-$.

\proclaim{Lemma 4.36} For any $2\le p<\infty$ and for $t$
sufficiently large, $|\int \tilde\mu\widetilde w - \int
\boldsymbol\mu \bold{w}| \le ct^{-\frac34+\frac1{2p}}$.
\endproclaim

\demo{Proof} Write $\int \tilde\mu \widetilde w - \int
{\boldsymbol\mu \bold w} = \int \widetilde w-\bold{w} +
\int(\boldsymbol\mu-I) (\widetilde w -\bold{w}) +
\int(\tilde\mu-\boldsymbol\mu)\widetilde w \equiv \text{I} +
\text{II} + \text{III}$. Again, after translation $z\to z_0+z$, we
obtain by (4.18), $|\text{I}| \le ct^{-3/4}$. Again by
triangularity, $\text{II} =
\int(C_{\bold{w}_\theta}\boldsymbol\mu)(\widetilde w-\bold{w}) =
\int(C^+\boldsymbol\mu \bold{w}^-_\theta)(\widetilde w^+_\theta -
\bold{w}_\theta{}^+) + (C^- \boldsymbol\mu \bold {w}_\theta{}^+)
(\widetilde w^-_\theta - \bold{w}_\theta{}^-) \equiv \text{II}_+ +
\text{II}_-$. By Cauchy, (4.33), and the $z\to z_0+z$ translation
of (4.21)--(4.22), $|\text{II}_+| = |\int (C^+ \boldsymbol\mu
\bold{w}_\theta{}^-) C^-(\widetilde w^+_\theta -
\bold{w}_\theta{}^+)|\le \|C^+ {\boldsymbol\mu \bold
w}_\theta{}^-\|_{L^2} \|C^-(\widetilde w^+_\theta -
\bold{w}_\theta{}^+)\|_{L^2}\le$ $ct^{-\frac14}t^{-\frac12} =
ct^{-\frac34}$.

The estimate for $\text{II}_-$ is similar and therefore $|\text{II}| \le ct^{-\frac34}$. We compute
$$
\tilde\mu - {\boldsymbol\mu} = (1-C_{\bold{w}_\theta})^{-1} C_{\tilde w-\bold{w}} \tilde\mu = C_{\tilde w-\bold{w}} \tilde \mu + C_{\bold{w}_\theta}\phi,
$$
where $\phi = (1-C_{\bold{w}_\theta})^{-1} C_{\tilde w-\bold{w}} \tilde\mu$. Thus by triangularity and Cauchy,
$$\align
\text{III} &= \int (C_{\tilde w-\bold{w}}\tilde\mu) \widetilde w + \int (C_{\bold{w}_\theta}\phi)\widetilde w\\
&= -\int (C^+\tilde\mu(\widetilde w^-_\theta - \bold{w}_\theta{}^-)) C^- \widetilde w^+_\theta + \int(C^-\tilde\mu(\widetilde w^+_\theta - \bold{w}_\theta{}^+)) C^+\widetilde w^-_\theta\\
&\quad - \int (C^+\phi\bold{w}_\theta{}^-) C^- \widetilde w^+_\theta + \int(C^- \phi \bold{w}_\theta{}^+)C^+ \widetilde w^-_\theta\\
&\equiv A_+ + A_- + B_+ + B_-.
\endalign$$
By Lemma 4.35 and (the translate of) (4.21), $|A_+| \le
\|C^+\tilde\mu(\bold{w}_\theta{}^- - \widetilde
w^-_\theta)\|_{L^2} \|C^- \bold{w}_\theta{}^+\|_{L^2} \le
ct^{-\frac34+\frac1{2p}}$, with a similar estimate for $A_-$.

Now consider $B_\pm$. Again by Lemma 4.35, $\|C_{\tilde
w-\bold{w}} \tilde\mu\|_{L^2} \le ct^{-\frac12+\frac1{2p}}$ and
hence $\|\phi\|_{L^2} \le ct^{-\frac12+\frac1{2p}}$. Together with
(4.21), this implies $\|B_\pm\|_{L^2} \le
ct^{-\frac34+\frac1{2p}}$, and hence $|\text{III}| \le
ct^{-\frac34+\frac1{2p}}$. Adding up the estimates for I, II and
III, the result follows.$\qquad\square$
\enddemo

By (4.3) and the lemma above, $q = \frac1{2\pi} \int (\tilde
\mu\widetilde w)_{12} = \frac1{2\pi} \int(\boldsymbol\mu
\bold{w})_{12} + O(t^{-\frac34+\frac1{2p}})$. By the results of \S
2 and the residue calculation at the end of \S 3, we have
$\int\boldsymbol\mu \bold{w} = \int_{\Sigma^e} \mu^0w^0$ where
$\mu^0 \in I + L^2(\Sigma^e)$ is the solution of
$(1-C_{w^0})\mu^0=I$. Let $\mu^\#\in I+L^2(\Sigma^e)$ be the
solution of $(1-C_{w^\#})\mu^\# = I$. Write
$$\align
\int_{\Sigma^e} \mu^\#w^\# - \int_{\Sigma^e} \mu^0w^0 &= \int_{\Sigma^e}(w^\#-w^0) + \int_{\Sigma^e} (\mu^\#-I) (w^\#-w^0) + \int_{\Sigma^e} (\mu^\#-\mu^0)w^0\\
&= \text{I} + \text{II} + \text{III}.
\endalign$$
By (2.63), $|\text{I}| \le \frac{c}{t^{3/4}}$. From the formula $\mu^\# - I = (1-C_{w^\#})^{-1} (C_{w^\#}I)$, we have by (2.65) and the explicit form of $v^\#$ in (2.61), $\|\mu^\#-I\|_{L^2} \le \|(1-C_{w^\#})^{-1}\|_{L^2} \|C^-(v^\#-I)\|_{L^2} \le \frac{c}{t^{1/4}}$. Again by (2.63) we conclude $|\text{II}| \le \frac{c}{t^{1/4}} \cdot \frac{c}{t^{\frac14+\frac14}} = \frac{c}{t^{3/4}}$. Finally $\mu^\# -\mu^0 = (1-C_{w^0})^{-1} C^-(\mu^\#(v^\#-v^0))$ and so
$$\align
\|\mu^\#-\mu^0\|_{L^2} &\le c(\|(\mu^\#-I)\|_{L^2} \|v^\#-v^0\|_{L^\infty} + \|v^\#-v^0\|_{L^2})\\
&\le c\left(\frac1{t^{1/4}} \frac{1}{t^{1/4}} + \frac{1}{t^{\frac14+\frac14}}\right) \le \frac{c}{t^{1/2}}.
\endalign$$
As $\|w^0\|_{L^2} \le \frac{c}{t^{1/4}}$, we conclude that $|\text{III}|\le c/t^{3/4}$. Hence as $t\to\infty$,
$$
\int_{\Sigma^e} \mu^\#w^\# - \int_{\Sigma^e} \mu^0w^0 = O\left(\frac1{t^{3/4}}\right).
$$
Assembling the above results we conclude that for any $p\ge 2$
$$
q(x,t) = \frac1{2\pi} \int_{\Sigma^e}(\mu^\#w^\#)_{12} +
O(t^{-3/4+1/(2p)}).
$$
But the normalized RHP $(\Sigma^e,v^\#)$ can be solved explicitly
in terms of parabolic cylinder functions as in \cite{DIZ},
\cite{DZ2} and we find $\frac1{2\pi} \int_{\Sigma^e}
(\mu^\#w^\#)_{12} = q_{\roman {as}}(x,t)$ where $q_{\roman
{as}}(x,t)$ is precisely the form $t^{-1/2}\alpha(z_0)
e^{ix^2/4t-i\nu(z_0)\log 2t}$ in (1.2). This completes the proof
of Theorem~1.10 with $q(t=0)\in H^{1,1}$.

Finally, we observe that the error term only depends on the
$H^{1,0}$ norm of $r = \cl R^{-1}(q(t=0))$. However, as mentioned
in the Introduction, equation (1.1) is not well-posed in
$H^{0,1}$. Indeed, for $q(t=0)=q_0\in H^{0,1}$, if $q(t)\in
H^{0,1}$ for $t>0$, then $r(t)=\cl R(q(t))\in H^{1,0}$. But for
$t>0, \,r(t)= e^{-iz^2t}r(t=0)$, which does not lie in $H^{1,0}$
for general $r(t=0)=\cl R(q_0),\, q_0\in H^{0,1}$. However,
standard contraction-mapping methods using the Strichartz-type
estimate
$$\|e^{-i\loz H_0}f\|_{L^\infty(dx)\otimes L^4_\rb(dt)}=\left(\int^
\infty_{-\infty}(\|u(t)\|_{L^{\infty}(dx)}^4\right)^{1/4}\le
\|f\|_{L^2(dx)},\tag 4.37$$ together with the standard estimate
$$\|e^{-itH_0}f\|_{L^q(dx)}\le (4\pi
|t|)^{1/q-1/2}\|f\|_{L^p(dx)},\quad 1/p+1/q=1,\quad 1\le p\le2\, ,
\tag 4.38$$ yields the following well-known result (see, for
example, [CW][GV][CSU] and the references therein).

\proclaim {Theorem 4.39} Let $q_0\in L^2(\rb)$ be given. Then
there exists a unique solution $q$ of (1.9) with $q\in
C(\rb_+,L^2(dx))\bigcap \,({L^\infty(dx)\otimes
L^4_{\roman{loc}}(dt)})$, and $q(t=0)=q_0$. Also
$\|q(t)\|_{L^2(dx)}$ is conserved.
\endproclaim

\demo{Remark}
Theorem 4.39 is commonly referred to as
well-posedness in $L^2$. This is, of course, something of a
misnomer as we need to place additional restrictions on the
solution in order to prove uniqueness. This nomenclature, however,
is well-established in the literature.
\enddemo
We show that our solution $q=q(x,t)$ obtained from the RHP (4.1)
via (4.3) solves (1.9) in the sense of Theorem 4.39. In view of
the preceding remarks (recall also Remark 3.28) this means that we
have in fact proved the following extended version of Theorem
1.10. \proclaim {Theorem 4.40}Let $q(t),\ t\ge0$, solve (1.9) with
$q_0=q(t=0)\in H^{0,1}\subset L^2$, as in Theorem 4.39. Fix
$0<\kappa<1/4$. Then as
$t\to\infty$,$$q(x,t)=t^{-1/2}\alpha(z_0)e^{ix^2/(4t)-i\nu(z_0)\log2t}+
O\left(t^{-(1/2+\kappa)}\right),$$ where $\alpha$ and $\nu$ are
given in terms of $r=\sr(q_0)$ as before. The error term $
O\left(t^{-(1/2+\kappa)}\right)$ is uniform for all
$x\in\rb$.\endproclaim

Note first that $\int \mu(w^+_\theta+w^-_\theta)=\int
(\mu-I)(w^+_\theta+w^-_\theta)+ \int
(w^+_\theta+w^-_\theta)$=I+II. By (2.30), we have, $|I|\le c
\frac{{\norm2r}^2}{1-\rho}.$ As the non-zero entries of $\int
(w^+_\theta+w^-_\theta)$ are proportional to $e^{-itH_0}\check r$
and its complete conjugate, we have the bound
$\|II\|_{L^\infty(dx)\otimes L^4_\rb(dt)}\le c\norm2r$ by (4.37).
This implies the bound
$$\|q\|_{L^\infty(dx)\otimes L^4_{[0,T]}(dt)}\le c \frac{\norm2r^2
T^{1/4}}{1-\rho} + c\norm2r\tag 4.41$$ for any $T>0$. Now if
$q_{n0}\in H^{1,1}$, say, n=1,2,..., let $q_n(t)\in
C(\rb_+,H^{1,1})$ be the global solution of (1.9) in $H^{1,1}$
described in the Introduction. For such solutions we have the
conservation law
$$\|q_n(\loz,t)\|_{L^2(dx)}=-\frac1{2\pi}\int \log(1-|r_n(z)|^2)dz
,\tag 4.42$$ where $r_n=\cl R(q_{n0})$. For $T>0$, define the norm
$$|||u|||_{[0,T]}\equiv \max_{0\le t\le T}||u(t)||_{L^2(dx)}+\|u\|_{L^\infty(dx)\otimes
 L^4_{[0,T]}(dt)}.$$ Then a standard calculation applied to (1.9)
 using (4.38)(4.39) implies $$ |||q_n-q_m|||_{[0,T]}\le c
 ||q_{0n}-q_{0m}||_{L^2(dx)} + c T^{1/2}(\|q_n\|_{L^\infty(dx)\otimes
 L^4_{[0,T]}(dt)}+\|q_m\|_{L^\infty(dx)\otimes
 L^4_{[0,T]}(dt)})^2|||q_n-q_m|||_{[0,T]}\tag 4.43$$ for
 $n,m\ge1$. If $q_0\in H^{0,1}$ is given, we an choose $q_{n0}\in
 H^{1,1}$ such that $q_{n0}\rightarrow q_0$ in $H^{0,1}$. But then
 $r_n=\cl R(q_{n0})\rightarrow r=\cl R(q_0)$ in $H^{1,0}$ and hence  ${\|q_n\|_{L^\infty(dx)\otimes
 L^4_{[0,T]}(dt)}}\le c $ for all $n\ge1$ by (4.41). It follows
then by (4.43) that the $q_n$'s converge in $|||\cdot|||_{[0,T]}$,
say $q_n\rightarrow q^0$, for some (sufficiently small) $T_0>0$ .
Taking the limit as $n\rightarrow\infty$ in (1.9), we see, in
particular, that $q^0$ solves (1.9) in the sense of Theorem 4.39,
at least for $0\le t \le T_0$. But then we can repeat the argument
starting at $t=T_0$, and we conclude that $q^0$ solves (1.9) for
all $t\ge0$. Finally, as $r_n\rightarrow r$ in $H^{1,0}$, we can
take the limit $n\rightarrow \infty$ in I+II to conclude that for
$t>0,\, \|q_n(t)-q(t)\|_{L^\infty(dx)}\rightarrow 0$. But for
$t>0,\,
\|q_n(t)-q^0(t)\|_{L^2(dx)}\le|||q_n-q^0|||_{[0,t]}\rightarrow 0$
and we conclude that $q(t)=q^0(t)$, and hence $q$ solves (1.9) in
the sense of Theorem 4.39. This concludes the proof of Theorem
4.40.

\demo{Remark}

As $q_n(t)\rightarrow q^0(t)=q(t)$ in $L^2(dx)$, we learn from
(4.42) that
$$\|q(\loz,t)\|_{L^2(dx)}=-\frac1{2\pi}\int \log(1-|r(z)|^2)dz=\roman {const.}$$
It is an interesting fact that we do not seem able to derive this
conservation law for $q_0\in H^{0,1}$ directly from the inverse
scattering formalism.
\enddemo

\end